\RequirePackage{fix-cm}
\documentclass[smallextended]{svjour3}       
\smartqed  
\usepackage{graphicx}
 \usepackage{epstopdf}
\usepackage{mathptmx}
\usepackage{latexsym}
\usepackage{float}
\usepackage{amssymb}
\usepackage{array}
\usepackage{multirow}
\usepackage{amsmath}
\usepackage{enumerate}
\usepackage{color}
\usepackage{bbding}
\usepackage{subfigure}
\usepackage{amscd}
\usepackage{color}
\usepackage{booktabs}
\usepackage{indentfirst}
\usepackage [latin1]{inputenc}
\allowdisplaybreaks
\numberwithin{equation}{section}
\def\proof{\par{\it Proof}. \ignorespaces}

\def\RR{\mathbb{R}}
\def\NN{\mathbb{N}}
\def\calP{\mathcal{P}}
\def\QED{\fbox{}}
\journalname{Numer Algo}
\begin{document}

\title{Hamiltonian Boundary Value Methods (HBVMs) for functional
differential equations with piecewise continuous arguments\thanks{This work is supported by NSFC (No. 12201537,12071403,12101525) and Research Support Funding of Xiangtan University (No. 22QDZ23,21QDZ16).}}
\titlerunning{}        
\author{Gianmarco Gurioli$^{\textbf{1}}$ \and Weijie Wang$^{\textbf{2}}$  \and Xiaoqiang Yan$^{\textbf{2}},^{\textbf{3}}$}
\institute{Gianmarco Gurioli  \\
           \email{gianmarco.gurioli@unifi.it}\\
           Weijie Wang \\
           \email{WeijieWang@smail.xtu.edu.cn}\\
            (\Envelope)~Xiaoqiang Yan  \\
           \email{xqyan1992@xtu.edu.cn} \\
  \at
           $^1$~ Universit\`{a} degli Studi di Firenze, Florence, Italy. Member of GNCS-INdAM.
  \at
           $^2$~ School of Mathematics and Computational Science $\&$ Hunan Key Laboratory for Computation and Simulation in Science and Engineering, Xiangtan University, Hunan 411105, China.
 \at
           $^3$~Key Laboratory of Intelligent Computing and Information Processing of Ministry of Education, Xiangtan University, Hunan 411105, China.
}
\date{Received: date / Accepted: date}

\maketitle

\begin{abstract}
In this paper, a class of high-order methods to numerically solve Functional Differential Equations with Piecewise Continuous Arguments (FDEPCAs) is discussed. The framework stems from the expansion of the vector field associated with the reference differential equation along the shifted and scaled Legendre polynomial orthonormal basis, working on a suitable extension of Hamiltonian Boundary Value Methods. Within the design of the methods, a proper generalization of the perturbation results coming from the field of ordinary differential equations is considered, with the aim of handling the case of FDEPCAs. The error analysis of the devised family of methods is performed, while a few numerical tests on Hamiltonian FDEPCAs are provided, to give evidence to the theoretical findings and show the effectiveness of the obtained resolution strategy.

\medskip


\end{abstract}

\section{Introduction}\label{section1}
In this paper, we consider initial value problems for Functional Differential Equations with Piecewise Continuous Arguments (FDEPCAs) of the form:
\begin{equation}\label{Pb}
\left\{\begin{array}{llll}
\dot{y}\left( t \right) = f\left( t,y\left( t \right),y(\lfloor t\rfloor) \right),~~~t \in \left[ {0,b} \right]\\[2\jot]
y\left( 0 \right) = {y_0}\\[2\jot]
\end{array}\right.,
\end{equation}
with ${y_0}$, $y \in {\mathbb{R}^m}$, $b$ a positive scalar, $f$: $\left[ {0,b} \right] \times {\mathbb{R}^m} \times {\mathbb{R}^m} \to {\mathbb{R}^m}$ a suitably regular function. As conventional, the over dot ``~$\dot{~}$~'' denotes the first-order derivative with respect to the time variable $t$ and $\lfloor .\rfloor$ the floor\footnote{i.e.,  for a given $t\in[0,b]$, $\lfloor t\rfloor$ is the greatest integer smaller than or equal to $t$.} function. As stated in \cite{Wiener93}, a solution of \eqref{Pb} is a function $y\left( t \right)$ fulfilling the following conditions:
\begin{itemize}
  \item[$\bullet$] ${y\left( t \right)}$ is continuous on $t \in \left[ {0,b} \right]$;
  \item[$\bullet$] the first-order derivative $\dot{y}\left( t \right)$ exists at each point $t \in \left[ {0,b} \right]$, with the possible exception of the point $\lfloor t \rfloor \in \left[ {0,b} \right]$, where one-sided derivatives exist;
  \item[$\bullet$] problem \eqref{Pb} is satisfied on each interval $\left[ {n,n + 1} \right)\subseteq \left[ {0,b} \right]$ ($n\in\mathbb{N}$) with integral endpoints.
\end{itemize}
 
It is worth noticing that the differential equation in \eqref{Pb} can be seen as a particular Delay Differential Equation (DDE, see \cite{BellZenn2003}), since the argument $\lfloor t\rfloor$ of the delay term $y(\lfloor t\rfloor)$ is such that $\lfloor t\rfloor \le t$, for all $t\in[0,b]$, yet with the substantial difference of presenting a discontinuous (piecewise continuous and piecewise constant) delay. FDEPCAs are a class of significant differential equation widely used to model phenomena arising from several scientific fields, such as biology (\cite{Akhmet06,Altintan06,Busenberg12,Gopalsamy92}), physics
 (\cite{Esedog06,Feng20,Sobolev64}), mechanics (\cite{Bereketoglu10,Chiu19,Esmaeilzadeh20,Kumar13,Driver63}), control science (\cite{Chiu19,Kolmanovskii12}) or population dynamics (\cite{Gopalsamy98}). 
 
 A substantial theory of FDEPCAs with a delay term of the form $y(t)=[t]$ ($[\cdot]$ denoting the integer part function) was approached in  \cite{Mysh1977} and more deeply addressed by the reference literature \cite{Cooke84,Cooke86,Gopalsamy98,Wiener93}, also considering differential equations with multiple piecewise continuous delay terms. Since for such equations a closed-form solution is not generally known, many authors have worked on the development of classical numerical schemes, mainly relying on Runge-Kutta methods (RKs) \cite{Liang14,Liu04,Liu09,Liu15,Lv11,Wang08,Wang13} and Linear Multi-step Methods (LMMs) \cite{Song10,Wen05}, studying the numerical stability, the preservation of oscillations and convergence properties. A third line of research has exploited methods in between LMMs and RKs: Boundary Value Methods (BVMs), dating back to \cite{Brugnano97-1,Brugnano98-1} within the range of Ordinary Differential Equations (ODEs), and Block BVMs (BBVMs) \cite{Brugnano97-2}. In this direction, the authors in \cite{Li17,Li19} applied BBVMs to solve \eqref{Pb} and neutral equations with piecewise constant arguments, while Zhang \& Yan \cite{Zhang21} further extended BBVMs to approximate the solutions of nonlinear delay differential equations with algebraic constraint and piecewise continuous arguments.
 
 More in general, BVMs have been applied to linear Hamiltonian systems of ODEs (see, e.g., \cite{Brugnano97-3}) to get symplectic integrators of arbitrary high order in the class of LMMs and then turned into the class of Hamiltonian Boundary Value Methods (HBVMs) to deal with the energy-preservation of Hamiltonian systems of ODEs, by building on the so-called discrete line integral approach \cite{Brugnano16-1}. The framework joins the methods relying on a local Fourier expansion of the vector field associated to the evolving system along an orthonormal basis, that specifically refers to the Legendre polynomial basis $\{P_j\}_{j\ge 0}$:
\begin{equation}\label{Leg}
P_j\in\Pi_j, \qquad \int_0^1P_i(x)P_j(x)\,dc =\delta_{ij},\qquad i,j\in\mathbb{N},
\end{equation}
with $\Pi_j$ the vector space of polynomials of degree $j$ and $\delta_{ij}$ the Kronecker delta.
In the past years, efficient energy-conserving HBVMs have been proposed to numerically solve a wide range of Hamiltonian problems, encountering ``drift-free'' methods for polynomial Hamiltonian systems \cite{Brugnano09-1}, energy-preserving resolutions of Poisson systems \cite{Brugnano12-1} and of Hamiltonian systems \cite{Brugnano 15-1}, also with holonomic constraints \cite{Brugnano18-2}, energy and quadratic invariants preserving strategies \cite{Brugnano12-4,Brugnano18-1} and multiple-invariants conserving schemes \cite{Brugnano14-1}, with an eye on their efficient implementation \cite{Brugnano16-1}. Recently, the papers \cite{Brugnano20-1,Brugnano18-4,Brugnano19-1} shown that HBVMs can be seen as spectral methods in time, highlighting their potential in the resolution of highly oscillatory problems and Hamiltonian partial differential equations \cite{Brugnano18-3,Brugnano18-5,Brugnano19-3,Brugnano19-4,Brugnano19-5}. Basing on the above mentioned background, the authors in \cite{Brugnano22-1} eventually employed HBVMs to solve constant delay differential equations with a continuous delay term, providing a perturbation analysis for such class of differential equations.

To the best of our acknowledge, there has not been any attempt to solve problem \eqref{Pb} using extensions of HBVMs, which in turn seem to be promising in the resolution of the class of delay differential equations in \cite{Brugnano22-1}. In view of this, we will work on such an extension, taking into account a crucial starting point: the derivation of a suitable perturbation analysis for the reference problem \eqref{Pb}. In fact, expanding the right-hand side term of the ODE in \eqref{Pb} along the Fourier expansion given by the polynomial basis \eqref{Leg}, up to a chosen number of terms, leads to a more numerically easy to solve differential problem, that constitutes a projection of the original differential problem onto a finite dimensional vector space. Consequently, the obtained problem can be seen as a perturbation of the original one and, hence, assessing the relationship between the solutions of the two problems becomes of paramount importance. For such a reason, suitable perturbation results for problem \eqref{Pb} are here theoretically provided.

To these aims, the structure of the paper is as follows: in Section \ref{section2} we state the reference problem, obtaining a polynomial approximation of the problem solution by means of the truncated Fourier expansion based on \eqref{Leg}, while the main perturbation result is proved in Section \ref{section3}. The subsequent Section \ref{section4} is devoted to the accuracy analysis of the obtained polynomial approximation to the solution. This analysis is fundamental to derive our extension of HBVMs to problem \eqref{Pb} and perform its error analysis, that is the object of Section \ref{section5}, together with hints about its efficient implementation. Section \ref{section6} then reports some numerical tests for problem \eqref{Pb} and, finally, concluding remarks and future perspectives are given in Section \ref{section7}.

\section{Statement of the problem}\label{section2}

Without loss of generality, to set the  numerical resolution of the problem, let us refer to the simpler form:
\begin{equation}\label{eq:ode1}
\left\{\begin{array}{llll}
\dot{y}\left( t \right) = f\left( y\left( t \right),y(\lfloor t\rfloor) \right),~~~t \in \left[ {0,b} \right]\\[2\jot]
y\left( 0 \right) = {y_0}\\[2\jot]
\end{array}\right.
\end{equation}
of problem \eqref{Pb}. Moreover, splitting the interval $[0,b]$ into $N\in\mathbb{N}$ subintervals with fixed mesh points
\begin{equation}\label{eq:mesh}
t_n = nh, \qquad n=0,\dots,N;  \qquad h = \nu^{-1}, \qquad 0\neq \nu\in\NN,
\end{equation}
we can restrict the solution of \eqref{eq:ode1} to each subinterval $[t_{n-1},t_n]$ ($n=1,2,...,N$), by formally setting

\begin{equation}\label{eq:hsign}
\hat\sigma_n(ch) := y(t_{n-1}+ch) \equiv \hat\sigma(t_{n-1}+ch), \qquad c\in[0,1],
\end{equation}
where the function $\hat\sigma(t)\equiv y(t)$ is introduced for usefull notational purposes. Consequently, one has:
\[
\hat\sigma_n(h)=y(t_n), \qquad n=0,1,\dots,N.
\]
Without loss of generality, it is hereafter assumed that $N=K\nu$, for some $K\in\mathbb{N}$.

Generalizing the arguments in \cite {Brugnano16-1}, we shall use the local Fourier expansion of the right-hand of the ODE in \eqref{eq:ode1} along the shifted and scaled Legendre polynomial orthogonal basis \eqref{Leg}, obtaining, for $n=1,...,N$:
\begin{equation}\label{eq:hsig1}
\left\{\begin{array}{llll}
&\dot{\hat\sigma}_n(ch) = f\Bigl(\hat\sigma_n(ch),\hat\sigma_{\lfloor \frac{n+c}\nu\rfloor\nu}(0)\Bigr)=\sum_{j\ge 0} P_j(c)\gamma_j(\hat\sigma_n,\hat\sigma_{\lfloor \frac{n+c}\nu\rfloor\nu}), ~~ c\in[0,1]\\&\hat\sigma_n(0) = y(t_{n-1})
\end{array}\right.,
\end{equation}
\noindent
with
\begin{equation}\label{eq:gammaj}
\gamma_j(z,w) := \int_0^1 P_j(\zeta)f(z(\zeta h),w(0))\, d\zeta,\quad j\ge 0,
\end{equation}
for any suitably regular functions $z,w:[0,h]\rightarrow\RR^m$. Clearly, according to \eqref{eq:mesh}--\eqref{eq:hsign},
\[
\hat\sigma_{\lfloor \frac{n+c}\nu\rfloor\nu}(0) = \hat\sigma\left(t_{\lfloor \frac{n-1+c}\nu\rfloor\nu}(0)\right)=\hat\sigma\left(\lfloor (n-1+c)/\nu\rfloor\right).
\]

Integrating both sides of the ODE in \eqref{eq:hsig1} with respect to $c\in[0,1]$, one obtains  the following formal expression for the solution of \eqref{eq:hsig1}:
\begin{equation}\label{eq:hsig}
\hat\sigma_n(ch) = y(t_{n-1})+h\sum_{j\ge0} \int_0^c P_j(x)\,dx\,\gamma_j(\hat\sigma_n,\hat\sigma_{\lfloor \frac{n+c}\nu\rfloor\nu}), \qquad c\in[0,1]
\end{equation}
and, by virtue of the orthonormality property in \eqref{Leg},
\begin{equation}\label{eq:ytn}
\hat\sigma_n(h) = y(t_n) = y(t_{n-1}) + h\gamma_0(\hat\sigma_n,\hat\sigma_{\lfloor \frac{n+1}\nu\rfloor\nu}) \equiv \hat\sigma(t_n).
\end{equation}

We now look for a piecewise polynomial approximation $\sigma(t)$, of the solution of (\ref{eq:ode1}), such that $\sigma_n\in\Pi_s$, for $n=1,\dots,N$,  obtained by truncating the local Fourier expansion in \eqref{eq:hsig1} to the first $s$ terms. Recalling \eqref{eq:gammaj}, problem \eqref{eq:hsig1} is indeed replaced by the following:

\begin{equation}\label{eq:hsig2}
\left\{\begin{array}{llll}
&\dot{\sigma}_n(ch) = \sum_{j=0}^{s-1} P_j(c)\gamma_j(\sigma_n,\sigma_{\lfloor \frac{n+c}\nu\rfloor\nu}), \qquad c\in[0,1]\\&\sigma_n(0) = \sigma(t_{n-1})
\end{array}\right.,
\end{equation}
where, as usual, we set
\begin{equation}\label{eq:sign}
\sigma_n(ch) \equiv \sigma(t_{n-1}+ch), \qquad c\in[0,1].
\end{equation}
Similarly to the derivation of \eqref{eq:hsig}, $\sigma_n$ can be formally written as:

\begin{equation}\label{eq:sig}
\sigma_n(ch) = \sigma_n(0)+ h\sum_{j=0}^{s-1} \int_0^c P_j(x)\,d x \, \gamma_j(\sigma_n,\sigma_{\lfloor \frac{n+c}\nu\rfloor\nu}), \qquad c\in[0,1],
\end{equation}
and, when $c=1$ is considered (compare with (\ref{eq:ytn})),
\begin{equation}\label{eq:yn-1}
\sigma_n(h) = \sigma_n(0) + h\gamma_0(\sigma_n,\sigma_{\lfloor \frac{n+1}\nu\rfloor\nu}) \equiv \sigma(t_n), 
\end{equation}
being an approximation of the solution $\hat\sigma(t)$ at the grid point $t_n$, $n=1,...,N$.
\begin{remark}\label{Rem:HBVMinf}
Equations \eqref{eq:sig}--\eqref{eq:yn-1} define a so-called HBVM($\infty$,$s$) method on the $n$-th time interval $[t_{n-1},t_n]$ ($n=1,...,N$), that can be recast as a continuous Runge-Kutta scheme of order $2s$ (see, e.g., \cite{Brugnano16-1}).
\end{remark}
\noindent
We end this section reporting few preliminary results related to \eqref{eq:gammaj}, to be considered for later use.

\begin{lemma}\label{thm:Gh} (\cite{Brugnano20-1}) Let $G:[0,h]\rightarrow V$, with $V$ a vector space, admit a Taylor expansion at 0. Then, for all $j\ge0$:
$$\int_0^1 P_j(\zeta) G(\zeta h)\,d\zeta = O(h^j).$$
\end{lemma}
\proof  By virtue of (\ref{Leg}):
\begin{eqnarray*}
\int_0^1 P_j(\zeta) G(\zeta h)\,d\zeta &=& \int_0^1 P_j(\zeta) \sum_{k\ge0} \frac{G^{(k)}(0)}{k!}(\zeta h)^k\,d\zeta =  \sum_{k\ge0} \frac{G^{(k)}(0)}{k!}h^k\int_0^1P_j(\zeta)\zeta^k\,d\zeta\\
&=& \sum_{k\ge j} \frac{G^{(k)}(0)}{k!}h^k\int_0^1P_j(\zeta)\zeta^k\,d\zeta = O(h^j),
\end{eqnarray*}
with $G^{(0)}=G$ and $G^{(k)}$, $k> 0$, denoting the derivative of $G$ of order $k$. ~~~\QED

\begin{corollary}\label{cor:gamj} With reference to (\ref{eq:gammaj}), one has: \,$\gamma_j(z,w) = O(h^j)$, for $j\ge 0$.
\end{corollary}
\proof
The thesis is a straight consequence of the previous Lemma \ref{thm:Gh}, setting $G(\zeta h):=f(z(\zeta h),w(0))$.~~~\QED

\begin{lemma}\label{lemma3}
Let ${G:\left[ {0,h} \right] \to V}$, with $V$ a vector space and ~$G\left( 0 \right) \ne 0$, admit a Taylor expansion at $0$. Then, for all $c  \in \left( {0,1} \right)$, and for all $j \ge 0$:
\[
\int_0^c {{P_j }\left( \zeta  \right)G\left( {\zeta  h} \right)} d\zeta  = O\left( {{h^0 }} \right)= O\left( 1 \right).
\]
\end{lemma}
\proof
\begin{equation*}
\begin{aligned}
\int_0^c {{P_j }} \left( \zeta \right)G\left( {\zeta  h} \right)d\zeta   &= \int_0^c {{P_j }} \left( \zeta \right)\sum\limits_{k \ge 0} {\frac{{{G^{\left( k \right)}}\left( 0 \right)}}{{k!}}} {\left( {\zeta h} \right)^k}d\zeta\\
 &= \sum\limits_{k \ge 0} {\frac{{{G^{\left( k \right)}}\left( 0 \right)}}{{k!}}} {h^k}\int_0^c  {{P_j }} \left( \zeta  \right){\zeta ^k}d\zeta
.\end{aligned}
\end{equation*}
Since, for $c\in(0,1)$, $\int_0^c {{P_j}} \left( \zeta  \right){\zeta ^k}d\zeta \ne 0$ and, by assumption, $G\left( 0 \right) \ne 0~$, it follows that:
$$\int_0^c {{P_j }} \left(\zeta  \right)G\left( {\zeta h} \right)d\zeta  = O\left( {{h^0}} \right) = O\left( 1 \right).~~~\QED$$

\section{Perturbation results}\label{section3}

In this section we establish perturbation results for problem \eqref{eq:ode1}, that will be used for the accuracy analysis of Section \ref{section4} and Section \ref{section5}. To this aim, without loss of generality, we shall discuss them on a local problem defined on the time interval $[\xi,T]$, with $\xi\in[0,1)$, $T\in(\xi,1)$:
\begin{equation}\label{eq:odexi}
\left\{\begin{array}{llll}
\dot y(t) = f(y(t),\phi), \quad t\in[\xi,T]\\[1\jot]
y\left( \xi \right) = \eta
\end{array}\right.,
\end{equation}
in which the constant memory term $y(\lfloor t \rfloor)$ is denoted by $\phi$. As one may easily notice, problem \eqref{eq:odexi} can be seen as a perturbation of the original problem \eqref{eq:ode1} in $[0,T]$, occurring for $\xi=0$ and $\eta = y_0\in\RR^m$. Let us denote the solution of \eqref{eq:odexi} by
\begin{equation}\label{eq:ysol}
y(t)\equiv y(t,\xi,\eta,\phi),
\end{equation}
to stress its dependance on the arguments $t$, the starting time $\xi$, the initial condition $\eta$ and the delay term $\phi$.
Moreover, let us make use of the notations
\begin{equation}\label{eq:F12}
F_1(z,w) = \frac{\partial}{\partial z} f(z,w), \qquad F_2(z,w) = \frac{\partial}{\partial w} f(z,w),
\end{equation}
to briefly address the partial derivatives of $f$.

The following theorem states perturbation results with respect to all the arguments $t$, $\xi$, $\eta$, $\phi$ of the solution \eqref{eq:ysol} (we shall hereafter refer to either the first or the second notation in \eqref{eq:ysol}, depending on the needs).

\begin{theorem}\label{thm:pertres}
With reference to the solution (\ref{eq:ysol}) of problem (\ref{eq:odexi}), for $t\in[\xi,T]$, one has:
\begin{align*}
&a)~\frac{\partial}{\partial t}y(t) =f(y(t),\phi),\\
&b)~\frac{\partial}{\partial \eta} y(t) = \Phi(t,\xi,\eta,\phi,y),\\
&c)~\frac{\partial}{\partial \xi} y(t) = -\Phi(t,\xi,\eta,\phi,y)f(\eta,\phi),
\end{align*}
where $\Phi(t,\xi,\eta,\phi,y)$ is the solution of the variational problem
\begin{equation}\label{eq:variational}
\left\{\begin{array}{llll}
\dot\Phi(t,\xi,\eta,\phi,y) = F_1(y(t),\phi)\Phi(t,\xi,\eta,\phi,y),\quad t\in[\xi,T]\\[1\jot]
\Phi(\xi,\xi,\eta,\phi,y)= I_m
\end{array}\right.,
\end{equation}
being $I_m\in\RR^{m\times m}$ the $(m\times m)$ unit matrix.\\
Further,
$$d)~\frac{\partial}{\partial \phi} y(t) = \int_{\xi}^t \Phi(t,s,\eta,\phi,y)F_2(y(s),\phi)\,d s.$$
\end{theorem}

\proof The statement a) is straightforward due to the problem definition in \eqref{eq:odexi}. Let us consider a perturbation $\delta\eta\in\RR^m$ on the initial condition $\eta$, leading to the perturbed initial condition $\tilde{\eta}=\eta+\delta\eta$. Problem \eqref{eq:odexi} thus results into the following perturbed one:
\begin{equation}\label{eq:odexi2}
\left\{\begin{array}{llll}
\dot{\tilde{y}}(t) = f(\tilde{y}(t),\phi), \quad t\in[\xi,T]\\[1\jot]
\tilde{y}\left( \xi \right) = \tilde\eta
\end{array}\right.,
\end{equation}
whose solution can be denoted by $\tilde y(t)= y(t,\xi,\tilde\eta,\phi)$.
Setting $z(t)=\tilde{y}(t)-y(t)$ and recalling the left-hand side equation in \eqref{eq:F12}, together with the ODEs in \eqref{eq:odexi} and \eqref{eq:odexi2}, the Taylor expansion then gives
\[
\dot z(t) = \dot{\tilde{y}}(t) - \dot{y}(t) = F_1(y(t),\phi) z(t)+ O\left(\|z(t)\|^2\right).
\]
Neglecting the higher order terms, we get the variational problem
\[
\left\{\begin{array}{llll}
\dot z(t) = F_1(y(t),\phi) z(t),\quad t\in[\xi,T]\\[1\jot]
z(\xi) = \delta\eta
\end{array}\right.,
\]
with solution
\begin{equation}\label{eq:z}
z(t)=y(t,\xi,\eta+\delta\eta,\phi)-y(t,\xi,\eta,\phi)=\Phi(t,\xi,\eta,\phi,y) \delta\eta,\quad t\in[\xi,T],
\end{equation}
such that
\[
\dot\Phi(t,\xi,\eta,\phi,y) = F_1(y(t),\phi)\Phi(t,\xi,\eta,\phi,y),\quad t\in[\xi,T],
\]
and
\[
z(\xi)=\Phi(\xi,\xi,\eta,\phi,y)\delta\eta=\delta\eta,
\]
implying that $\Phi(\xi,\xi,\eta,\phi,y)=I_m$. The proof of the point b) is then completed noticing that, by \eqref{eq:z}:
\[
\frac{\partial}{\partial \eta} y(t) =\Phi(t,\xi,\eta,\phi,y).
\]
Concerning the point c), let us consider a generic $t^*$ in the interval $(\xi,T]$ and set $y^*=y(t^*,\xi,\eta,\phi)$; consequently,
~$\eta = y(\xi,t^*,y(t^*,\xi,\eta,\phi),\phi)$~ and, therefore, the identity
$$y^* = y(t^*,\xi,y(\xi,t^*,y^*,\phi),\phi)$$ holds true.
By considering the derivative of both members with respect to $\xi$, recalling the thesis of the statement b) and denoting by $0_m$ the $m$-dimensional null vector, we have:
\begin{eqnarray*}
0_m&=&\frac{d}{d\xi} y^* =\frac{\partial}{\partial\xi} y(t^*,\xi,\eta,\phi) + \left.\frac{\partial}{\partial\eta}y(t^*,\xi,\eta,\phi)\frac{\partial}{\partial t}y(t,\xi,\eta,\phi)\right|_{t=\xi}\\
&=&\frac{\partial}{\partial\xi} y(t^*) + \Phi(t^*,\xi,\eta,\phi,y)f(\eta,\phi).
\end{eqnarray*}
The statement eventually follows because of the arbitrariness of $t^*\in(\xi,T]$.\\
Regarding the last point d), due to \eqref{eq:odexi}, \eqref{eq:F12} and the thesis of the statement b), one has:
\[
\begin{split}
\frac{d}{d t}\frac{\partial}{\partial \phi} y(t)&=\frac{\partial}{\partial \phi}\dot{y}(t) =\frac{\partial}{\partial \phi} f(y(t),\phi) \\
&=\frac{\partial}{\partial z}f(z,\phi)\Big|_{z=y(t)}\frac{\partial}{\partial \phi} y(t)+\frac{\partial}{\partial w}f(y(t),w)\Big|_{w=\phi}\frac{\partial}{\partial \phi} \phi\\
&= F_1(y(t),\phi)\frac{\partial}{\partial \phi} y(t) + F_2(y(t),\phi).
\end{split}
\]
Further,
\[
\frac{\partial}{\partial \phi}y(\xi) = \frac{\partial}{\partial \phi}\eta = O_m,
\]
where $O_m$ stands for the $(m\times m)$ zero matrix. For the sake of simplicity, let us denote $\frac{\partial}{\partial \phi} y(t)$ by $y_{\phi}(t)$. We then search for a solution of the problem
\begin{equation}\label{eq:ic}
\left\{\begin{array}{llll}
\dot{y}_{\phi}(t)= F_1(y(t),\phi) y_{\phi}(t)+F_2(y(t),\phi),\quad t\in[\xi,T]\\[1\jot]
y_{\phi}(\xi) = O_m
\end{array}\right.,
\end{equation}
of the form
\begin{equation}\label{eq:yform}
y_{\phi}(t)=W(t)C(t),
\end{equation}
with $C(t)\equiv C(t,\xi,\eta,\phi,y)\in\RR^{m\times m}$ and a nonsingular matrix $W(t)\equiv W(t,\xi,\eta,\phi,y)\in\RR^{m\times m}$, for $t\in[\xi,T]$, satisfying
\begin{equation}\label{eq:W}
\dot{W}(t)=F_1(y(t),\phi) W(t),\quad t\in[\xi,T].
\end{equation}
Setting
\begin{equation}\label{eq:Psi}
\Psi(t,\xi)\equiv \Psi(t,\xi,\eta,\phi,y) =W(t)W^{-1}(\xi),\quad t\in[\xi,T],
\end{equation}
we highlight that  \eqref{eq:Psi} and \eqref{eq:W} imply
\[
\dot{\Psi}(t,\xi)=\dot{W}(t)W^{-1}(\xi)=F_1(y(t),\phi)\Psi(t,\xi),\quad t\in[\xi,T],
\]
with $\Psi(\xi,\xi)=W(\xi)W^{-1}(\xi)=I_m$ and, hence, $\Psi=\Phi$ in \eqref{eq:variational}.
Moreover, considering the time derivative of both sides of \eqref{eq:yform} and recalling \eqref{eq:ic}--\eqref{eq:W}, we have that
\begin{equation}\label{eq:doty}
\begin{split}
\dot{y}_{\phi}(t) &= \dot{W}(t)C(t)+W(t)\dot{C}(t) = F_1(y(t),\phi) W(t)C(t)+W(t)\dot{C}(t)\\
& =F_1(y(t),\phi)y_{\phi}(t)+ W(t)\dot{C}(t),\quad t\in[\xi,T],
\end{split}
\end{equation}
while the fulfillment of the initial condition $y_{\phi}(\xi)=W(\xi)C(\xi)$ gives
\begin{equation}\label{eq:C0}
C(\xi)=W^{-1}(\xi)y_{\phi}(\xi).
\end{equation}
As a consequence, imposing that $W(t)\dot{C}(t)=F_2(y(t),\phi)$ in \eqref{eq:doty}, for $t\in[\xi,T]$, we obtain:
\[
C(t)=C(\xi)+\int_{\xi}^t W^{-1}(s)F_2(y(s),\phi)\,ds,\quad t\in[\xi,T].
\]
Plugging the above equation into \eqref{eq:yform} and recalling \eqref{eq:C0}, the definition of $\Psi=\Phi$ in \eqref{eq:Psi} and the initial condition in \eqref{eq:ic}, finally reads:
\[
\begin{split}
y_{\phi}(t) &= W(t)C(t)\\
&=W(t)W^{-1}(\xi)y_{\phi}(\xi)+\int_{\xi}^t W(t)W^{-1}(s)F_2(y(s),\phi)\,ds\\
& =\Phi(t,\xi,\eta,\phi,y)y_{\phi}(\xi)+\int_{\xi}^t  \Phi(t,s,\eta,\phi,y)F_2(y(s),\phi)\,ds\\
& = \int_{\xi}^t \Phi(t,s,\eta,\phi,y)F_2(y(s),\phi)\,ds, \quad t\in[\xi,T],
\end{split}
\]
that concludes the proof.~~~\QED\\

\noindent
A significant corollary is stated below, referring to any chosen vector norm $|\cdot |$.
\begin{corollary}\label{cor:pert} With reference to the solution (\ref{eq:ysol}) of problem \eqref{eq:odexi}, for any $\delta\eta$, $\delta\phi\in\RR^m$, one has:
\[
\begin{split}
y(t,\xi,\eta+\delta\eta,\phi+\delta\phi)=&~ y(t,\xi,\eta,\phi) + \Phi(t,\xi,\eta,\phi+\delta\phi,y)\delta\eta \\
&+ \int_\xi^t\Phi(t,s,\eta,\phi,y)F_2(y(s),\phi)\,d s\,\delta\phi,\\
&+(t-\xi)O(|\delta\eta|^2 +|\delta\phi|^2), \qquad t\in[\xi,T].
\end{split}
\]
\end{corollary}
\proof
By virtue of the statements b) and d) of Theorem \ref{thm:pertres}, one has:
\[
\begin{split}
y(t,\xi,\eta+\delta\eta,\phi+\delta\phi)&= y(t,\xi,\eta+\delta\eta,\phi+\delta\phi)\pm y(t,\xi,\eta,\phi+\delta\phi)\pm y(t,\xi,\eta,\phi)\\
& = [y(t,\xi,\eta+\delta\eta,\phi+\delta\phi)-y(t,\xi,\eta,\phi+\delta\phi)]\\
&~~~~+[y(t,\xi,\eta,\phi+\delta\phi)- y(t,\xi,\eta,\phi)]+y(t,\xi,\eta,\phi)\\
& =\frac{\partial y}{\partial \eta}(t,\xi,\eta,\phi+\delta\phi)\delta\eta+\frac{\partial y}{\partial \phi}(t,\xi,\eta,\phi)\delta\phi\\
&~~~~+y(t,\xi,\eta,\phi)+(t-\xi)O(|\delta\eta|^2+|\delta\phi|^2)\\
& = y(t,\xi,\eta,\phi) + \Phi(t,\xi,\eta,\phi,y)\delta\eta+\\
&~~~~+\int_\xi^t\Phi(t,s,\eta,\phi,y)F_2(y(s),\phi)\,d s\,\delta\phi\\
&~~~~+(t-\xi)O(|\delta\eta|^2+|\delta\phi|^2), \qquad t\in[\xi,T].~~~\QED
\end{split}
\]
The following relevant auxiliary result is also proved.

\begin{lemma}\label{lem:fih}
With reference to the fundamental matrix function $\Phi$ introduced in Theorem~\ref{thm:pertres}, for any sufficiently small $h>0$, one has:
\[
\Phi(\xi+h,\xi,\eta,\phi,y) = I_m+O(h), \qquad \int_\xi^{\xi+h}\Phi(\xi+h,s,\eta,\phi,y)F_2(y(s),\phi)\,d s = O(h).
\]
\end{lemma}
\proof From \eqref{eq:variational} and allowing $\Phi$ for Taylor expansion at $\xi$, we have that
\[
\begin{split}
\Phi(\xi+h,\xi,\eta,\phi,y) &= \Phi(\xi,\xi,\eta,\phi,y)+h\dot{\Phi}(t,\xi,\eta,\phi,y)\Big|_{t=\xi}+O(h^2)\\
& = I_m+hF_1(y(\xi),\phi)\Phi(\xi,\xi,\eta,\phi,y)+O(h^2)\\
& = I_m+hF_1(y(\xi),\phi)+O(h^2)=I_m+O(h),\\
\end{split}
\]
so that the left hand side equality of the statement is proved. Setting $s=\xi+ch$, $c\in[0,1]$, then gives:
\[
\begin{split}
&\int_\xi^{\xi+h}\Phi(\xi+h,s,\eta,\phi,y)F_2(y(s),\phi)\,ds\\
&=h\int_0^1\Phi(\xi+h,\xi+ch,\eta,\phi,y)F_2(y(\xi+ch),\phi)\,dc=h\cdot O(h^0)=O(h),
\end{split}
\]
where the penultimate equality is due to Lemma \ref{thm:Gh}, so the statement is fully proved.~~~\QED

\section{Accuracy analysis}\label{section4}
We are now in the position to discuss the accuracy of the polynomial approximation $\sigma(t)$ at the grid-points (\ref{eq:mesh}), moving from the following preliminary result.
\begin{lemma}\label{thm:lem0}
With reference to \eqref{eq:mesh}--\eqref{eq:yn-1} and \eqref{eq:ysol}, assume that $\hat\sigma_{\lfloor \frac{n}\nu\rfloor\nu}(0)=\sigma_{\lfloor \frac{n}\nu\rfloor\nu}(0)=\sigma(\lfloor (n-1)/\nu\rfloor)=\phi$, then:
\[
y(t_n,t_{n-1},\sigma(t_{n-1}),\phi)-y(t_n,t_n,\sigma(t_n),\phi)=O(h^{2s+1}),
\]
\[
y(t_{n-1}+ch,t_{n-1},\sigma(t_{n-1}),\phi)-y(t_{n-1}+ch,t_{n-1}+ch,\sigma(t_{n-1}+ch),\phi)=O(h^{s+1}),\quad c\in(0,1).
\]
\end{lemma}

\proof Recalling \eqref{eq:hsig2}, \eqref{eq:sign}, \eqref{eq:yn-1} and the theses of Theorem \ref{thm:pertres}, Lemma \ref{thm:Gh}, Lemma \ref{lemma3} and Corollary \ref{cor:gamj}, one has that, for $c\in(0,1]$:
\[
\begin{split}
&y(t_{n-1}+ch,t_{n-1},\sigma(t_{n-1}),\phi)-y(t_{n-1}+ch,t_{n-1}+ch,\sigma(t_{n-1}+ch),\phi) =\\
&=y(t_{n-1}+ch,t_{n-1},\sigma_n(0),\phi)-y(t_{n-1}+ch,t_{n-1}+ch,\sigma_n(ch),\phi) \\
&=-\int_0^{ch}\frac{d}{d t}y(t_{n-1}+ch,t_{n-1}+t,\sigma_n(t),\phi)\,d t\\
&=-\int_0^{ch}\Bigl( \left.\frac{\partial}{\partial \xi} y(t_{n-1}+ch,\xi,\sigma_n(t),\phi)\right|_{\xi=t_{n-1}+t}+
\left.\frac{\partial}{\partial \eta} y(t_{n-1}+ch,t_{n-1}+t,\eta,\phi)\right|_{\eta=\sigma_n(t)}\dot\sigma_n(t)\Bigr)\,d t\\
&=\int_0^{ch} \Phi(t_{n-1}+ch,t_{n-1}+t,\sigma_n(t),\phi,y)\Bigl[f(\sigma_n(t),\phi)-\dot\sigma_n(t)\Bigr]\,d t\\
&=h\int_0^c \Phi(t_{n-1}+ch,t_{n-1}+\tau h,\sigma_n(\tau h),\phi,y)\Bigl[f(\sigma_n(\tau h),\phi)-\dot\sigma_n(\tau h)\Bigr]\,d \tau\\
&=h\int_0^c \underbrace{\Phi(t_{n-1}+ch,t_{n-1}+\tau h,\sigma_n(\tau h),\phi,y)}_{=:\,G(\tau h)}\left[\sum_{j\ge 0}P_j(\tau) \gamma_j(\sigma_n,\phi)-\sum_{j=0}^{s-1}P_j(\tau) \gamma_j(\sigma_n,\phi)\right]\,d \tau\\
&=h\int_0^c G(\tau h)\sum_{j\ge s}P_j(\tau ) \gamma_j(\sigma_n,\phi)\,d \tau=
h\sum_{j\ge s} \int_0^c P_j(\tau)G(\tau h)\,d \tau\,\gamma_j(\sigma_n,\phi)\\
& = \left\{\begin{array}{llll}
 O(h^{s+1}),~~\mbox{if} ~c\in(0,1)\\[2\jot]
 O(h^{2s+1}),~\mbox{if}~ c=1
 \end{array}\right..~~~\QED
\end{split}
\]

The following theorem can thus be proved.
\begin{theorem}\label{thm:acc0}
With reference to \eqref{eq:mesh}--\eqref{eq:yn-1} and \eqref{eq:ysol}, for $n=1,...,\nu$, one has
\[
y(t_n)-\sigma(t_n)=y(t_{n-1})+\sigma(t_{n-1})+O(h^{2s+1}),
\]
\[
y(t_{n-1}+ch)-\sigma(t_{n-1}+ch)=\hat\sigma_n(ch)-\sigma_n(ch)=O(h^{s+1}),\quad c\in(0,1).
\]
\end{theorem}
\proof The proof is done via generalized induction on $n$, taking into account that $y(t_0)=\hat\sigma_1(0)=\sigma_1(0)=\sigma(t_0)=y_0=\phi$.\\
First, the statement holds true when $n=1$ is considered. In fact, for $c\in(0,1]$,
\[
\begin{split}
&y(t_{0}+ch)-\sigma(t_{0}+ch)=\hat\sigma_1(ch)-\sigma_1(ch)=y(ch,0, \hat\sigma_1(0),\phi)-y(ch,ch,\sigma_1(ch),\phi)\\
&=y(ch,0, \hat\sigma_1(0),\phi)-y(ch,0,\sigma_1(0),\phi)+y(ch,0,\sigma_1(0),\phi)-y(ch,ch,\sigma_1(ch),\phi)\\
&=y(ch,0,\sigma_1(0),\phi)-y(ch,ch,\sigma_1(ch),\phi)=\left\{\begin{array}{llll}
 O(h^{s+1}),~~\mbox{if} ~c\in(0,1)\\[2\jot]
 O(h^{2s+1}),~\mbox{if}~ c=1
 \end{array}\right.,
\end{split}
\]
where the last equality is due to Lemma \ref{thm:lem0}. Therefore, for $c=1$,
\[
y(t_1)-\sigma(t_1)=y(t_0)-\sigma(t_0)+O(h^{2s+1})=O(h^{2s+1}).
\]
Suppose now that the statement holds until a generic $n-1$ ($n\in\{2,...,\nu-1\}$) and let us prove that it also holds for $n$. One has:
\[
\begin{split}
&y(t_{n-1}+ch)-\sigma(t_{n-1}+ch)=\hat\sigma_n(ch)-\sigma_n(ch)\\
&=y(t_{n-1}+ch,t_{n-1}, \hat\sigma_n(0),\phi)-y(t_{n-1}+ch,t_{n-1}+ch,\sigma_n(ch),\phi)\\
&=\underbrace{y(t_{n-1}+ch,t_{n-1}, \hat\sigma_n(0),\phi)-y(t_{n-1}+ch,t_{n-1},\sigma_n(0),\phi)}_{=:E_{n,1}(ch)}\\
&~~~~\underbrace{+y(t_{n-1}+ch,t_{n-1},\sigma_n(0),\phi)-y(t_{n-1}+ch,t_{n-1}+ch,\sigma_n(ch),\phi)}_{=:E_{n,2}(ch)},
\end{split}
\]
where, from Lemma \ref{thm:lem0}: $E_{n,2}(ch)=\left\{\begin{array}{llll}
 O(h^{s+1}),~~\mbox{if} ~c\in(0,1)\\[2\jot]
 O(h^{2s+1}),~\mbox{if}~ c=1
 \end{array}\right..$\\
Moreover, from Corollary \ref{cor:pert} and Lemma \ref{lem:fih}:
\[
\begin{split}
E_{n,1}(ch)&=\Phi(t_{n-1}+ch,t_{n-1},\sigma_n(0),\phi,y)(\hat\sigma_n(0)-\sigma_n(0))+hO(|\hat\sigma_n(0)-\sigma_n(0)|^2)\\
&=\left(I_m+O(h)\right)(\hat\sigma_n(0)-\sigma_n(0))+hO(|\hat\sigma_n(0)-\sigma_n(0)|^2)\\
&=\hat\sigma_n(0)-\sigma_n(0)+O(h)(\hat\sigma_n(0)-\sigma_n(0))+hO(|\hat\sigma_n(0)-\sigma_n(0)|^2).
\end{split}
\]
Recalling that, due to the induction hypotesis,
\[
\hat\sigma_n(0)-\sigma_n(0)=y(t_{n-1})-\sigma(t_{n-1})=(n-1)O(h^{2s+1}),
\]
one has that,
for $c\in(0,1)$,
\[
y(t_{n-1}+ch)-\sigma(t_{n-1}+ch)=E_{n,1}(ch)+E_{n,2}(ch)=O(h^{s+1}),
\]
while, for $c=1$:
\[
y(t_n)-\sigma(t_n)=E_{n,1}(h)+E_{n,2}(h)=y(t_{n-1})-\sigma(t_{n-1})+O(h^{2s+1}).~~~\QED
\]
We can finally generalize the previous theorem to the case of $\nu<n\le N$.
\begin{theorem}\label{thm:acc1}
With reference to \eqref{eq:mesh}--\eqref{eq:yn-1} and \eqref{eq:ysol}, for $n=1,...,N$, one has
\[
y(t_n)-\sigma(t_n)=y(t_{n-1})+\sigma(t_{n-1})+O(h^{2s+1}),
\]
\[
y(t_{n-1}+ch)-\sigma(t_{n-1}+ch)=\hat\sigma_n(ch)-\sigma_n(ch)=O(h^{s+1}),\quad c\in(0,1).
\]
\end{theorem}
\proof The proof is done by generalized induction on $n$. It is worth noticing that for $n=1,...,\nu$, the statement holds by virtue of Theorem \ref{thm:acc0}. Consequently, let us assume that the thesis is true until $n=\kappa\nu$ ($\kappa\in\{1,...,K-1\}$) and let us prove the statement for $n=\kappa\nu+1,...,(\kappa+1)\nu$. We first notice that, for $n=\kappa\nu+1$,...,$(\kappa+1)\nu$, the memory terms are given by $\hat{\phi}_{\kappa\nu}:=y(t_{\kappa\nu})=y_{\kappa\nu}(h)$ and $\phi_{\kappa\nu}:=\sigma(t_{\kappa\nu})=\sigma_{\kappa\nu}(h)$.  We have:
\[
\begin{split}
&y(t_{n-1}+ch)-\sigma(t_{n-1}+ch)=\hat\sigma_n(ch)-\sigma_n(ch)\\
&=y(t_{n-1}+ch,t_{n-1}, \hat\sigma_n(0),\hat\phi_{\kappa\nu})-y(t_{n-1}+ch,t_{n-1}+ch,\sigma_n(ch),\phi_{\kappa\nu})\\
&=\underbrace{y(t_{n-1}+ch,t_{n-1}, \hat\sigma_n(0),\hat\phi_{\kappa\nu})-y(t_{n-1}+ch,t_{n-1},\sigma_n(0),\phi_{\kappa\nu})}_{=:E_{n,1}^{(\kappa)}(ch)}\\
&~~~~\underbrace{+y(t_{n-1}+ch,t_{n-1},\sigma_n(0),\phi_{\kappa\nu})-y(t_{n-1}+ch,t_{n-1}+ch,\sigma_n(ch),\phi_{\kappa\nu})}_{=:E_{n,2}^{(\kappa)}(ch)},
\end{split}
\]
for which Lemma \ref{thm:lem0} still provides that $E_{n,2}^{(\kappa)}(ch)=\left\{\begin{array}{llll}
 O(h^{s+1}),~~\mbox{if} ~c\in(0,1)\\[2\jot]
 O(h^{2s+1}),~\mbox{if}~ c=1
 \end{array}\right..$\\
Concerning $E_{n,1}^{(\kappa)}(ch)$, from Corollary \ref{cor:pert} and Lemma \ref{lem:fih} it follows that
\[
\begin{split}
E_{n,1}^{(\kappa)}(ch)&=\Phi(t_{n-1}+ch,t_{n-1},\sigma_n(0),\hat\phi_{\kappa\nu},y)(\hat\sigma_n(0)-\sigma_n(0))\\&
~~~~+\int_{t_{n-1}}^{t_{n-1}+ch}\Phi(t_{n-1}+ch,s,\sigma_n(0),\phi_{\kappa\nu},y)\,ds(\hat\phi_{\kappa\nu}-\phi_{\kappa\nu})\\
&~~~~+hO(|\hat\sigma_n(0)-\sigma_n(0)|^2+|\hat\phi_{\kappa\nu}-\phi_{\kappa\nu}|^2)\\
&=(\hat\sigma_n(0)-\sigma_n(0))+O(h)(\hat\sigma_n(0)-\sigma_n(0))+O(h)(\hat\phi_{\kappa\nu}-\phi_{\kappa\nu})\\
&~~~~+hO(|\hat\sigma_n(0)-\sigma_n(0)|^2+|\hat\phi_{\kappa\nu}-\phi_{\kappa\nu}|^2).
\end{split}
\]
Reminding that $n=\kappa\nu+1,...,n=(\kappa+1)\nu$, we have that $n=O(h^{-1})$ and $\nu=h^{-1}$, so that the induction hypothesis reads:
$$\hat\sigma_n(0)-\sigma_n(0)=y(t_{n-1})-\sigma(t_{n-1})=(n-1)O(h^{2s+1})=O(h^{2s}),$$ $$\hat\phi_{\kappa\nu}-\phi_{\kappa\nu}=\kappa\nu O(h^{2s+1})=O(h^{2s})$$
and, hence,
\[E_{n,1}^{(\kappa)}(ch)=
\left\{\begin{array}{llll}
 O(h^{2s}),~~\mbox{if} ~c\in(0,1)\\[2\jot]
 y(t_{n-1})-\sigma(t_{n-1})+O(h^{2s+1}),~\mbox{if}~ c=1
 \end{array}\right..
 \]
Consequently, for $c\in(0,1)$,
\[
y(t_{n-1}+ch)-\sigma(t_{n-1}+ch)=E_{n,1}^{(\kappa)}(ch)+E_{n,2}^{(\kappa)}(ch)=O(h^{s+1}),
\]
while, for $c=1$:
\[
y(t_n)-\sigma(t_n)=E_{n,1}^{(\kappa)}(h)+E_{n,2}^{(\kappa)}(h)=y(t_{n-1})-\sigma(t_{n-1})+O(h^{2s+1}).~~~\QED
\]

\section{Formulation of the method and error analysis}\label{section5}
As preannounced in the Remark \ref{Rem:HBVMinf}, the procedure defined by \eqref{eq:sig}--\eqref{eq:yn-1} does not provide a computable numerical method, since the integrals in the definition \eqref{eq:gammaj} of the Fourier coefficients in \eqref{eq:hsig2} and \eqref{eq:sig} have to be properly approximated. This can be done by considering a suitable quadrature formula, leading to a possibly different piecewise polynomial approximation $u(t)\in\Pi_{s}$ of the solution of \eqref{eq:ode1}, satisfying
\begin{equation}\label{eq:Pbu}
\left\{\begin{array}{llll}
&\dot{u}_n(ch) = \sum_{j=0}^{s-1} P_j(c)\overline\gamma_j(u_n,u_{\lfloor \frac{n+c}\nu\rfloor\nu}), \qquad c\in[0,1]\\&u_n(0) = u(t_{n-1})
\end{array}\right.,
\end{equation}
setting, for $c\in[0,1]$, $u_n(ch) \equiv u(t_{n-1}+ch)$ and with
\begin{equation}\label{eq:bargamma}
\begin{split}
\overline\gamma_j(u_n,u_{\lfloor \frac{n+c}\nu\rfloor\nu})&=\sum_{i=1}^k b_iP_j(c_i)f\Bigl(u_n(c_ih),u_{\lfloor(n+c_i)/\nu\rfloor\nu}(0)\Bigr)\\
&=\gamma_j(u_n,u_{\lfloor \frac{n+c}\nu\rfloor\nu})-\Delta_j(h),
\end{split}
\end{equation}
where $(c_i,b_i)$ are the abscissae and the weights of the quadrature formula, the exact integral $\gamma_j(u_n,u_{\lfloor \frac{n+c}\nu\rfloor\nu})$ is defined as in \eqref{eq:gammaj} and $\Delta_j(h)$ denotes the quadrature error.
Therefore,
\begin{equation}\label{eq:unch}
u_n(ch) = u_n(0)+ h\sum_{j=0}^{s-1} \int_0^c P_j(x)\,d x \, \overline\gamma_j(u_n,u_{\lfloor \frac{n+c}\nu\rfloor\nu}), \qquad c\in[0,1],
\end{equation}
providing, when $c=1$, the approximation of the solution $\hat\sigma(t)$ at the grid point $t_n$ given by
\begin{equation}\label{eq:unh}
u_n(h) = u_n(0) + h\overline\gamma_0(u_n,u_{\lfloor \frac{n+1}\nu\rfloor\nu}) \equiv u(t_n),\quad n=1,...,N.
\end{equation}
Regarding the quadrature error in \eqref{eq:bargamma}, the following result is reported.
\begin{theorem}\label{thm:quad}
(\cite{Brugnano20-1}) If the quadrature $(c_i,b_i)$, $i=1,...,k$ has order $q$ (so that the quadratura is exact for polynomial integrands up to degree $q-1$), then
\[
\Delta_j(h)=O(h^{q-j}),\qquad j=0,1,...,s-1.
\]
\end{theorem}
\noindent
Hereafter, let us consider the Gauss-Legendre quadrature formula, having abscissae $c_i$, $i=1,...,k$, placed at the zeros of the $P_k$ polynomial in \eqref{Leg} and weights $b_i$ such that
\begin{equation}\label{eq:GL}
P_k(c_i)=0,\qquad b_i=\int_0^1 \ell_i(x)\,dx,\qquad \ell_i(x)=\prod_{j\neq i} \frac{x-c_i}{c_i-c_j},\qquad i\in\{1,...,k\},
\end{equation}
that has order $q=2k$ (see, e.g., \cite{Brugnano16-1}).
Actually, when the above quadrature formula \eqref{eq:GL} is employed, \eqref{eq:bargamma}--\eqref{eq:unh} define the $n$-th step of an HBVM($k$,$s$) method, with $1\le s\le k$, that can be seen as the $n$-th integration step of a Runge-Kutta method, with stages
\begin{equation}\label{eq:Yin}
Y_i^n:=u_n(c_ih)\in\RR^m,\quad i=1,...,k.
\end{equation}
This feature can be shown evaluating \eqref{eq:unch} at $c_1,...,c_k$ and taking into consideration \eqref{eq:bargamma}, that leads to
\begin{equation}\label{eq:Yin2}
Y_i^n=u_n(0) + h\sum_{\ell=1}^k b_{\ell}\sum_{j=0}^{s-1} \int_0^{c_i} P_j(x)\,d x \,P_j(c_{\ell})f_{\nu}(Y_{\ell}^n),\quad i=1,...,k,
\end{equation}
where we set $f_{\nu}(u_n(c_ih)):=f\Bigl(u_n(c_ih),u_{\lfloor(n+c_i)/\nu\rfloor\nu}(0)\Bigr)$, since the memory term $u_{\lfloor(n+c_i)/\nu\rfloor\nu}(0)$ is known when performing the $n$-th integration step. Moreover, due to the orthonormality property in \eqref{Leg},
\[
u_n(h)=u_n(0)+h\sum_{i=1}^k b_if_{\nu}(Y_i^n).
\]
We have hence obtained a $k$ stage Runge-Kutta method with abscissae and weights $(c_i,b_i)$ ($i=1,...,k$) and Butcher \textit{tableau}
\begin{equation}\label{tableau}
\begin{array}{c|c}
c
&  \mathcal{I}_s\calP_s^\top\Omega \in\mathbb{R}^{k\times k} \\
\hline
 & b^{\top}\\
\end{array},
\end{equation}
setting
\begin{equation}\label{def:Omega}
c=\left(c_1,...,c_k\right)^\top\in\RR^k,\quad b=\left(b_1,...,b_k\right)^\top\in\RR^k,\quad \Omega=\begin{pmatrix}
b_1 & &  \\
 & \ddots &  \\
 & & b_k \\
\end{pmatrix}\in\RR^{k\times k};
\end{equation}
\noindent
$\calP_s\in\RR^{k\times s}$ and $\mathcal{I}_s\in\RR^{k\times s}$ the matrices with $(i,j)$ components given by
\begin{equation}\label{def:PsIs}
(\calP_s)_{ij}=P_{j-1}(c_i),\quad (\mathcal{I}_s)_{ij}=\int_0^{c_i}P_{j-1}(x)\,dx,\quad i\in\{1,...,k\},~~j\in\{1,...,s\}.
\end{equation}

We can now extend the results of Theorem \ref{thm:acc0} and Theorem \ref{thm:acc1} to study the accuracy of the approximations \eqref{eq:unch}--\eqref{eq:unh}. To start with, let us generalize Lemma \ref{thm:lem0}.

\begin{lemma}\label{thm:lem0.1}
With reference to \eqref{eq:mesh}--\eqref{eq:ytn}, \eqref{eq:Pbu}--\eqref{eq:unh} and \eqref{eq:ysol}, assume that $\hat\sigma_{\lfloor \frac{n}\nu\rfloor\nu}(0)=u_{\lfloor \frac{n}\nu\rfloor\nu}(0)=u(\lfloor (n-1)/\nu\rfloor)=\phi$ and that the quadrature formula in \eqref{eq:bargamma} has order $q\ge 2s$, then:
\[
y(t_n,t_{n-1},u(t_{n-1}),\phi)-y(t_n,t_n,u(t_n),\phi)=O(h^{2s+1}),
\]
\[
y(t_{n-1}+ch,t_{n-1},u(t_{n-1}),\phi)-y(t_{n-1}+ch,t_{n-1}+ch,u(t_{n-1}+ch),\phi)=O(h^{s+1}),\quad c\in(0,1).
\]
\end{lemma}

\proof Due to \eqref{eq:Pbu}, \eqref{eq:unch}, \eqref{eq:unh} and the theses of Theorem \ref{thm:pertres}, one has that, for $c\in(0,1]$:
\[
\begin{split}
&y(t_{n-1}+ch,t_{n-1},u(t_{n-1}),\phi)-y(t_{n-1}+ch,t_{n-1}+ch,u(t_{n-1}+ch),\phi) =\\
&=y(t_{n-1}+ch,t_{n-1},u_n(0),\phi)-y(t_{n-1}+ch,t_{n-1}+ch,u_n(ch),\phi) \\
&=-\int_0^{ch}\frac{d}{d t}y(t_{n-1}+ch,t_{n-1}+t,u_n(t),\phi)\,d t\\
&=-\int_0^{ch}\Bigl( \left.\frac{\partial}{\partial \xi} y(t_{n-1}+ch,\xi,u_n(t),\phi)\right|_{\xi=t_{n-1}+t}+
\left.\frac{\partial}{\partial \eta} y(t_{n-1}+ch,t_{n-1}+t,\eta,\phi)\right|_{\eta=u_n(t)}\dot u_n(t)\Bigr)\,d t\\
&=\int_0^{ch} \Phi(t_{n-1}+ch,t_{n-1}+t,u_n(t),\phi,y)\Bigl[f(u_n(t),\phi)-\dot u_n(t)\Bigr]\,d t\\
&=h\int_0^c \Phi(t_{n-1}+ch,t_{n-1}+\tau h,u_n(\tau h),\phi,y)\Bigl[f(u_n(\tau h),\phi)-\dot u_n(\tau h)\Bigr]\,d \tau\\
&=h\int_0^c \underbrace{\Phi(t_{n-1}+ch,t_{n-1}+\tau h,u_n(\tau h),\phi,y)}_{=:\,G(\tau h)}\left[\sum_{j\ge 0}P_j(\tau) \gamma_j(u_n,\phi)-\sum_{j=0}^{s-1}P_j(\tau) \overline{\gamma}_j(u_n,\phi)\right]\,d \tau.
\end{split}
\]
Moreover, from \eqref{eq:bargamma} and Theorem \ref{thm:quad}, together with Lemma \ref{thm:Gh}, Lemma \ref{lemma3}, Corollary \ref{cor:gamj} and the assumption $q\ge 2s$, it follows:
\[
\begin{split}
&y(t_{n-1}+ch,t_{n-1},u(t_{n-1}),\phi)-y(t_{n-1}+ch,t_{n-1}+ch,u(t_{n-1}+ch),\phi) =\\
&=h\int_0^c G(\tau h)\left[\sum_{j\ge 0}P_j(\tau) \gamma_j(u_n,\phi)-\sum_{j=0}^{s-1}P_j(\tau) \left(\gamma_j(u_n,\phi)-\Delta_j(h)\right)\right]\,d \tau\\
&=h\sum_{j\ge s}\int_0^cP_j(\tau ) G(\tau h) \,d \tau\underbrace{\gamma_j(u_n,\phi)}_{=O(h^j)}+h\sum_{j=0}^{s-1}\int_0^c G(\tau h) P_j(\tau ) \,d \tau\underbrace{\Delta_j(h)}_{=O(h^{q-j})}\\
& = \left\{\begin{array}{llll}
 O(h^{s+1}),~~\mbox{if} ~c\in(0,1)\\[2\jot]
 O(h^{2s+1}),~\mbox{if}~ c=1
 \end{array}\right..~~~\QED
 \end{split}
\]
\noindent
As a consequence, it is straightforward to prove the following result concerning the accuracy of the obtained HBVM($k$,$s$) class of methods.
\begin{theorem}\label{thm:acc3}
With reference to \eqref{eq:mesh}--\eqref{eq:ytn}, \eqref{eq:Pbu}--\eqref{eq:unh}, \eqref{eq:ysol} and assume that the quadrature formula in \eqref{eq:bargamma} has order $q\ge 2s$, for $n=1,...,N$, one has
\[
y(t_n)-u(t_n)=y(t_{n-1})+u(t_{n-1})+O(h^{2s+1}),
\]
\[
y(t_{n-1}+ch)-u(t_{n-1}+ch)=\hat\sigma_n(ch)-u_n(ch)=O(h^{s+1}),\quad c\in(0,1).
\]
\end{theorem}
\proof The proof parallels the one of Theorem \ref{thm:acc0} and Theorem\ref{thm:acc1}, formally replacing $\sigma$ by $u$ and making use of Lemma \ref{thm:lem0.1}.~~~\QED
\noindent
\begin{remark} Theorem \ref{thm:acc3} implies that the derived HBVM($k$,$s$), $1\le s\le k$, class of methods has order $2s$; in fact, recalling that $N=K\nu=K/h$, we have that $y(t_N)-u(t_N)=NO(h^{2s+1})=O(h^{-1})O(h^{2s+1})=O(h^{2s})$. In addition, it is worth noticing that, when $k=s$ is considered, HBVM($k$,$s$)=HBVM($s$,$s$) represents the $s$-stage Gauss-Legendre collocation method of order $2s$.
\end{remark}
We conclude this section highlighting that HBVMs can be regarded as spectral methods in time, when a sufficiently large value of $k$ is adopted (see \cite{Brugnano20-1,Brugnano19-1,Brugnano18-4,Brugnano19-5}). In more detail, let us recall the expansion \eqref{eq:hsig1}, with the definition \eqref{eq:gammaj} and the notation $\hat\gamma_j^n:=\gamma_j(\hat\sigma_n,\hat\sigma_{\lfloor \frac{n+c}\nu\rfloor\nu})$. Assuming $\dot{\hat{\sigma}}_n\in L^2([0,h])$, it follows that
\[
\|\dot{\hat{\sigma}}_n\|_{L^2}^2=\sum_{j\ge 0}\|\hat\gamma_j^n\|^2<\infty\quad \Rightarrow \quad \|\hat\gamma_j^n\|\rightarrow 0, ~\mbox{as}~ j\rightarrow\infty.
\]
Consequently, under a finite precision arithmetic with machine epsilon $\epsilon$ and choosing the degree $s$ of the polynomial approximation $\sigma_n$ in \eqref{eq:hsig2} such that
\begin{equation}\label{spectral1}
\forall j>s:~~\| \hat\gamma_j^n\|<\epsilon \cdot \max_{i=0,...,s-1}\|\hat\gamma_i^n\|,
\end{equation}
it follows that
\begin{equation}\label{spectral2}
\sigma_n(ch)\doteq \hat{\sigma}_n(ch),\quad c\in[0,1],
\end{equation}
where $\doteq$ stands for ``equal within roundoff errors''. In so doing, using the notations $\gamma_j^n:=\gamma_j(\sigma_n,\sigma_{\lfloor \frac{n+c}\nu\rfloor\nu})$ and $\overline{\gamma}_j^n:=\overline{\gamma}_j(u_n,u_{\lfloor \frac{n+c}\nu\rfloor\nu})$ (recall \eqref{eq:hsig2} and \eqref{eq:Pbu}), the $k>s$ and $s$ values have to be large enough to get $\gamma_j^n\doteq \overline{\gamma}_j^n$, for $j=0,...,s-1$, and \eqref{spectral1} holds true, so that \eqref{spectral2} holds as well. Suitable choices for $k>s$ have been investigated in \cite{Brugnano20-1}. A common option is $k=\max\{20,s+2\}$ (see \cite{Brugnano18-4}).

Apparently, the use of large values of $k$ seems to computationally mismatch with \eqref{eq:Yin}, since $k$ is the number of stages of the Runge-Kutta scheme \eqref{tableau} that have to be computed to obtain the approximation $u_n(h)$ of $y(t_n)$, for $n=1,...,N$. Nevertheless, for the sake of completeness, we remind that the underlying discrete problem \eqref{eq:Yin2} to be solved can be recast so that it has block dimension $s$, independently of $k$. To this scope, let us set $e=\left(1,...,1\right)^\top\in\RR^k$ the unit $k$-dimensional vector, and let us denote by
\[
Y^n:=\begin{pmatrix}
Y_1^n\\
\vdots\\
Y_k^n
\end{pmatrix}\in\RR^{km}
\]
the stage block vector having block dimension $k$ (recall \eqref{eq:Yin}). From \eqref{eq:Yin2} and recalling the definitions of $\Omega\in\RR^{k\times k}$, $\calP_s\in\RR^{k\times s}$ and $\mathcal{I}_s\in\RR^{k\times s}$ in \eqref{def:Omega}--\eqref{def:PsIs}, the following vectorial stage equation is provided:
\begin{equation}\label{step1}
Y^n=e\otimes u_n(0)+h\mathcal{I}_s\calP_s^\top\Omega\otimes I_mf_{\nu}(Y^n)\in\RR^{km},
\end{equation}
where $f_{\nu}(Y^n)=(f_{\nu}(Y_1^n),...,f_{\nu}(Y_k^n))^\top\in\RR^{km}$ and ``$\otimes$'' denotes the tensor product. Moreover, from \eqref{eq:bargamma}, we have that (see, \cite{Brugnano16-1})
\begin{equation}\label{step2}
\overline{\gamma}^n:=\begin{pmatrix}\overline{\gamma}_0^n\\ \vdots \\ \overline{\gamma}_{s-1}^n \end{pmatrix}=\calP_s^\top\Omega\otimes I_mf_{\nu}(Y^n)\in\RR^{sm},
\end{equation}
that is a useful form for the block vector with the $s$ Fourier coefficients vectors in \eqref{eq:unch} (recall, also, \eqref{eq:bargamma}). Consequently, \eqref{step1} can be recast as
\begin{equation}\label{step3}
Y^n=e\otimes u_n(0)+h\mathcal{I}_s\otimes I_m\overline{\gamma}^n\in\RR^{km}.
\end{equation}
Plugging \eqref{step3} into \eqref{step2} we finally get
\begin{equation}\label{eq:finaldiscrete}
\overline{\gamma}^n=\calP_s^\top\Omega\otimes I_mf_{\nu}(e\otimes u_n(0)+h\mathcal{I}_s\otimes I_m\overline{\gamma}^n)\in\RR^{sm},
\end{equation}
that provides a discrete problem equivalent to \eqref{step1}, with block dimension $s$, independently of $k$. Of course, once such equation is solved, the approximation $u_n(h)$ of $y(t_n)$ is given by (recall \eqref{eq:unh}):
\[
u_n(h)=u_{n}(0)+h\overline{\gamma}_0^n,\quad n=1,...,N.
\]
We refer to \cite{Brugnano16-1} for efficient nonlinear iterations for solving \eqref{eq:finaldiscrete}, including the common used blended iteration of \cite{Brugnano2002}.


\section{Numerical experiments}\label{section6}

In this section, we report a few numerical tests for Hamiltonian FDEPCAs, in support of the correctness and effectiveness of the theoretical findings. It is worth mentioning that the use the energy-preserving HBVMs has been extremely satisfactory, when dealing with Hamiltonian systems, as shown by a series of paper (see, e.g., \cite{Brugnano22-2,Brugnano16-1,Brugnano18-4,Brugnano18-6,Brugnano19-1}). More recently (see \cite{Brugnano22-1}), also their application to Hamiltonian DDE with continuous arguments has resulted to be beneficial. We thus focus on their employment in Hamiltonian FDEPCAs, showing that the presented HBVM framework can gain prominent advantages as well, also when comparing with the Gauss Legendre collocation method of the same order. To this scope, we consider the following class of Hamiltonian FDEPCAs, ispired by the popular issue of searching for period orbits of DDEs (see, e.g., the work of the authors in \cite{DosBar1986,KapYor1974,MallNuss2011,Nuss1973,Nuss1979,Wal1975}):
\begin{equation}\label{5.1}
\left(\dot{q}(t),\dot{p}(t)\right)^\top = J\otimes I_m \nabla \left(H(\dot{q}(t),\dot{p}(t))+\alpha H(q(\lfloor t\rfloor),p(\lfloor t\rfloor))\right),\quad t\in(0,T],
\end{equation}
 defined\footnote{We refer to the notation: $\nabla H(q,p)=(\frac{\partial}{\partial q}H(q,p),\frac{\partial}{\partial p} H(q,p))^\top \in\mathbb{R}^{2m}$.} by an Hamiltonian function:
\begin{equation}\label{5.2}
	H:(q,p)\in \mathbb{R}^{m}\times \mathbb{R}^{m}\rightarrow \mathbb{R},
\end{equation}
and the orthogonal skew-symmetric matrix $J=\begin{pmatrix} 0 & 1\\-1 & 0\end{pmatrix}$, where $I_m$ stands for the $(m\times m)$ identity matrix, $\alpha\in\mathbb{R}$, $0<T\in\mathbb{R}$. Equation \eqref{5.1} is equipped with the given initial conditions 
\begin{equation}\label{ICtest}
q_{0}=q(0)\in\mathbb{R}^m,\qquad p_{0}=p(0)\in\mathbb{R}^m.
\end{equation}
For all the instances of \eqref{5.1}--\eqref{ICtest} covered in the following subsections, the corresponding behavior of the Hamiltonian and of the solution with respect to time have been inferred applying the spectrally accurate higher order method HBVM(22,20) with a small step size. As a matter of fact, we show that the theoretical order of convergence of the employed HBVM methods is numerically confirmed and that their use allows to effectively reproduce the geometric features of the solution, as the numerical approximation of the Hamiltonian becomes accurate. At this regard, the use of $k$ values larger than $s$ in the applied HBVM($k$,$s$) method, thus distinguishing from the Gauss Legendre method of the same order $2s$, turns to be relevant. For a given step size value $h$, to estimate the last point error ${\varepsilon _N}(h)$ and the order $p$ of convergence, in what follows we refer to the following formulas:
\begin{equation}\label{5.3}
{\varepsilon _N}(h) = {\left\|u_N(h/2)-u_N(h) \right\|_\infty },\;\;\;\;\;\;p = {\log _2}\frac{{{\varepsilon _N}\left( h \right)}}{{{\varepsilon _N}\left( {\frac{h}{2}} \right)}},
\end{equation}
where $u_N(h/2)$ and $u_N(h)$ represent the numerical solution values at the last time instant using the step size $h/2$ and $h$, respectively, while $\|\cdot\|_{\infty}$ denotes the infinity vector norm.

The presented HBVM methods have been implemented in Matlab (R2023b) and run on a 2.2 GHz Intel Core i7 dual-core computer with 8 GB of memory.

\subsection{Problem 1}\label{NT_Pb1}

Consider
	\begin{equation}\label{5.4}
	     \begin{array}{lllll}
			m=1,~~~~~~~~~~~H(q,p)=\frac{1}{4}(q^{4}+p^{4}),~~\\
			 \alpha=10^{-2},~~~~~~q_{0}=\sqrt{2},~~~~~~ p_{0}=0,
		\end{array}
	\end{equation}
to give the first example based on \eqref{5.1}--\eqref{ICtest}.

To give numerical evidence of the accuracy result of Theorem \ref{thm:acc3}, Table \ref{table1} reports the estimates of the last point error and of the order $p$ of convergence when using HBVM$(2,2)$, HBVM$(10,2)$ and HBVM($15$,$3$), to solve the problem \eqref{5.4} with step sizes $h = \frac{1}{{10q}}$ ($q = 2$, $4$, $8$, $16$). Clearly, all the methods exhibit the right order $2s$ of convergence and, for each value of $h$, the error values computed by HBVM$(2,2)$ are larger than the ones obtained by HBVM$(10,2)$.

\begin{table}[H]
		\centering
		\caption{Estimates of the last point error $\varepsilon_N(h)$ and of the convergence order $p$ of HBVM$(2,2)$, HBVM$(10,2)$ and HBVM($15$,$3$) to problem \eqref{5.1}--\eqref{ICtest} with \eqref{5.4}, referring to the time interval $[0,2]$ and step sizes $h= \frac{1}{{10q}}$ ($q = 2$, $4$, $8$, $16$).}
\begin{tabular}{lllllll}
\hline
\multicolumn{1}{c}{\multirow{2}{*}{$h$}} & \multicolumn{2}{c}{HBVM($2$,$2$)}                      & \multicolumn{2}{c}{HBVM($10$,$2$)}                      & \multicolumn{2}{c}{HBVM($15$,$3$)}                      \\ \cline{2-7} 
\multicolumn{1}{c}{}                   & \multicolumn{1}{c}{$\varepsilon_N(h)$} & \multicolumn{1}{c}{$p$} & \multicolumn{1}{c}{$\varepsilon_N(h)$} & \multicolumn{1}{c}{$p$} & \multicolumn{1}{c}{$\varepsilon_N(h)$} & \multicolumn{1}{c}{$p$} \\ \hline
$1/20$                                   & 5.1981e-06                   & \multicolumn{1}{c}{---} & 2.7882e-06                   & \multicolumn{1}{c}{---} & 4.5716e-09                   & \multicolumn{1}{c}{---} \\ \hline
$1/40$                                   & 3.2665e-07                   & 3.9922e+00              & 1.7570e-07                   & 3.9882e+00              & 7.1969e-11                   & 5.9892e+00              \\ \hline
$1/80$                                   & 2.0443e-08                   & 3.9980e+00              & 1.1004e-08                   & 3.9970e+00              & 1.1251e-12                   & 5.9993e+00              \\ \hline
$1/160$                                  & 1.2781e-09                   & 3.9995e+00              & 6.8809e-10                   & 3.9992e+00              & 1.9984e-14                   & 5.8151e+00              \\ \hline
\end{tabular}\label{table1}
\end{table}

\begin{figure}[!h]
\begin{center}
\includegraphics[width=.45\textwidth]{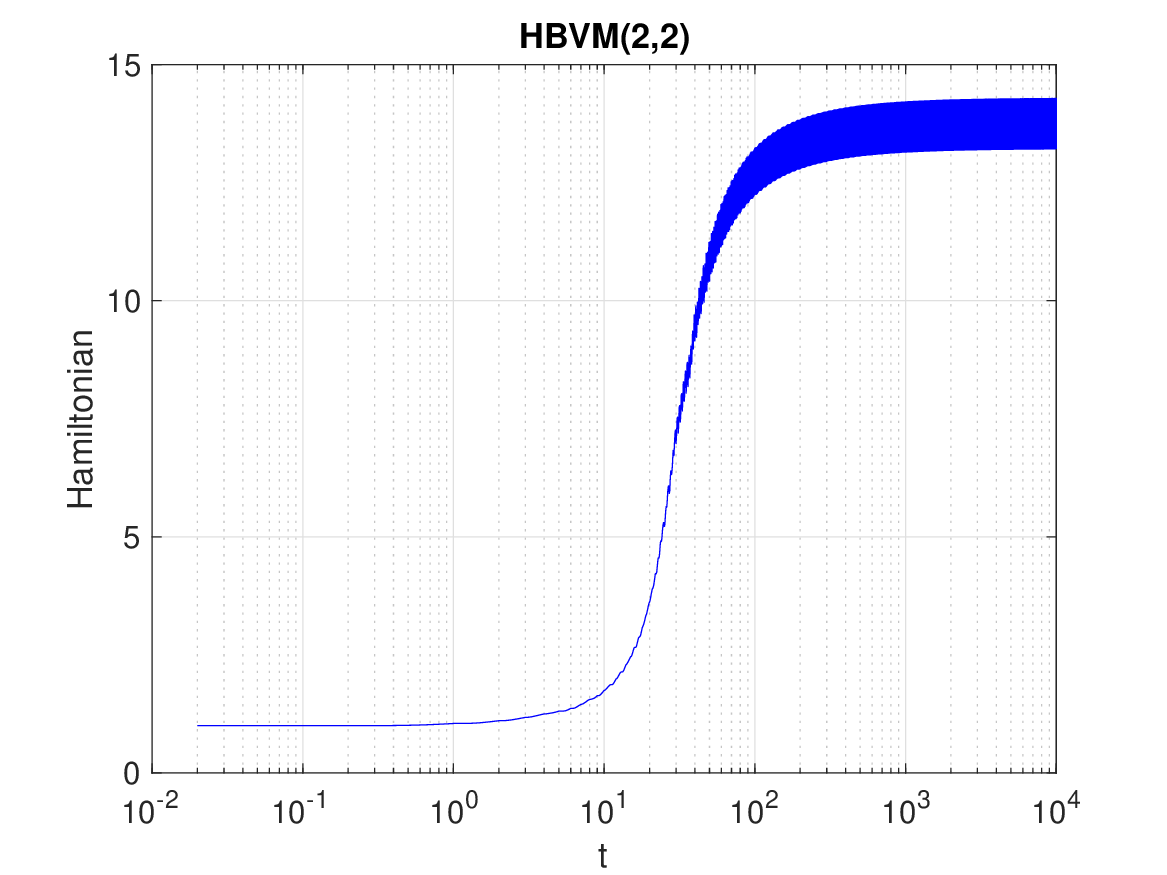}
\includegraphics[width=.45\textwidth]{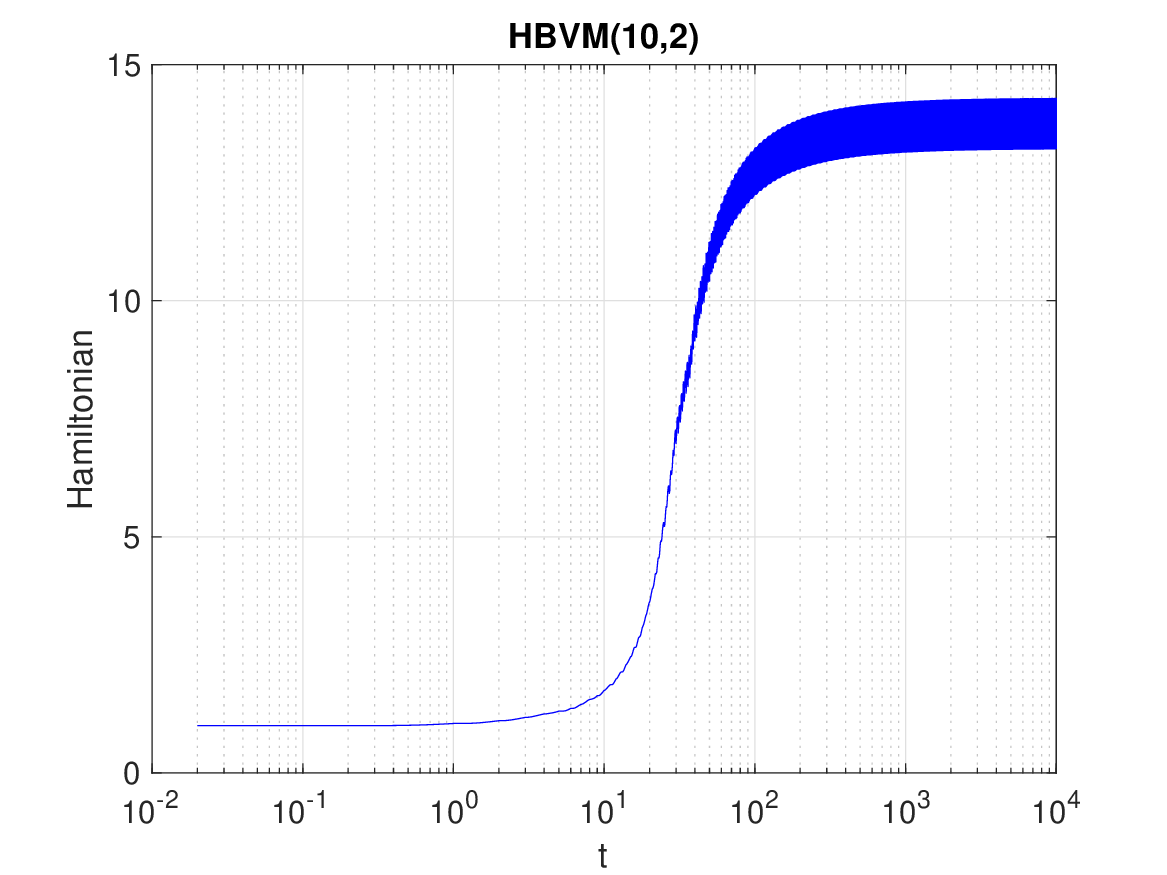}
\includegraphics[width=.45\textwidth]{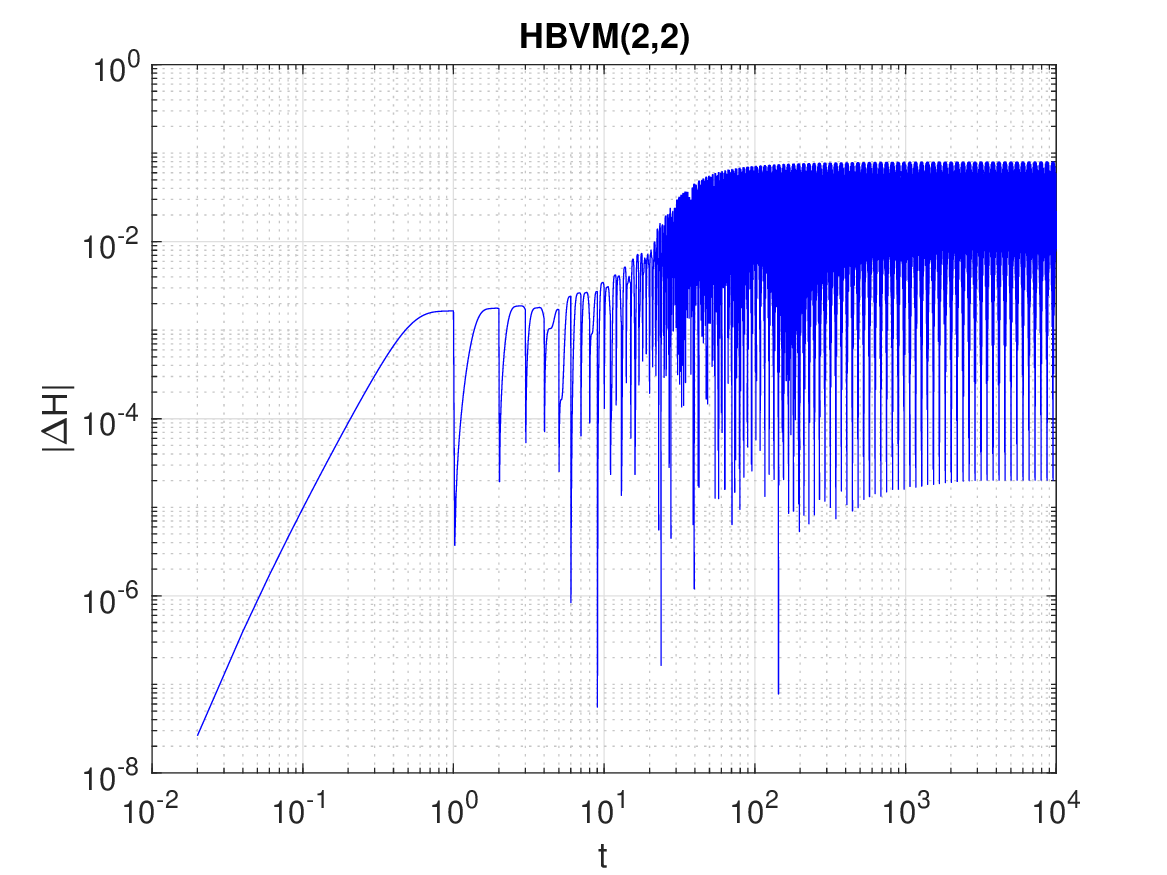}
\includegraphics[width=.45\textwidth]{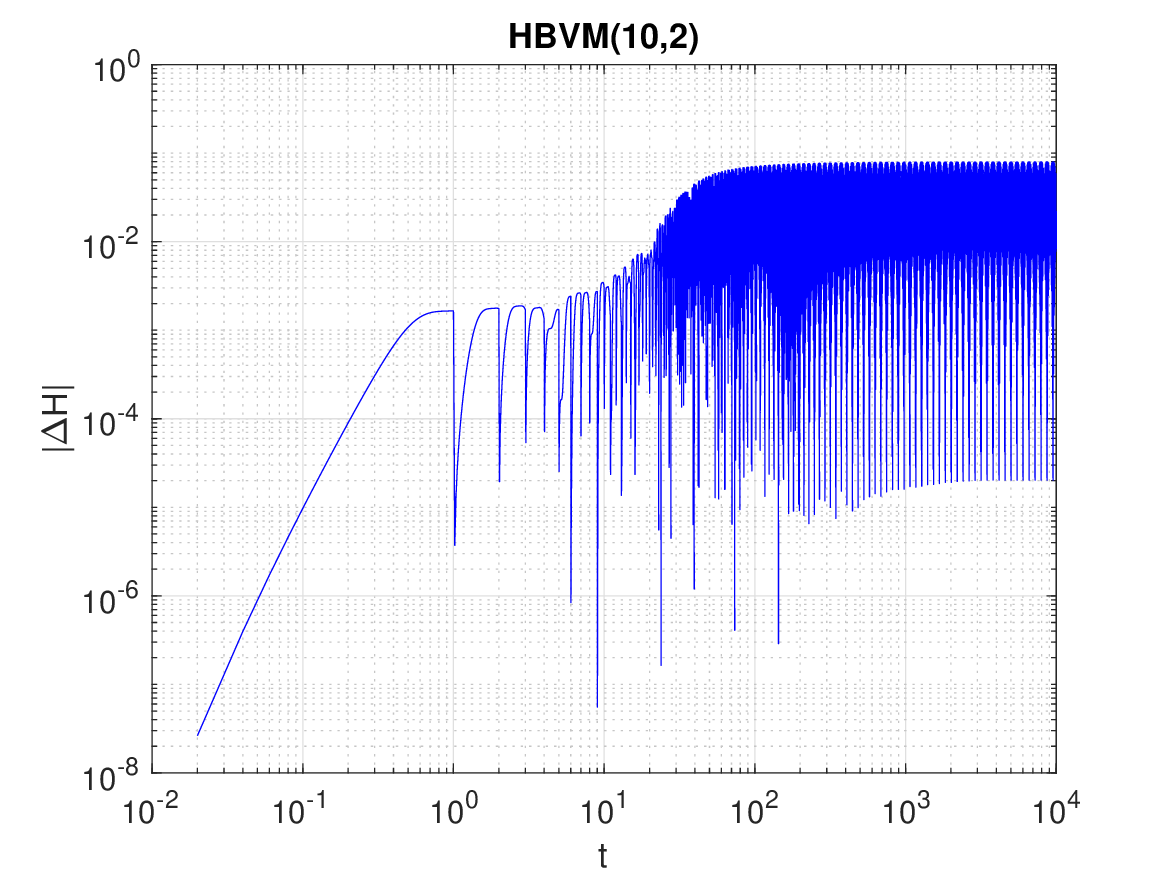}
\includegraphics[width=.45\textwidth]{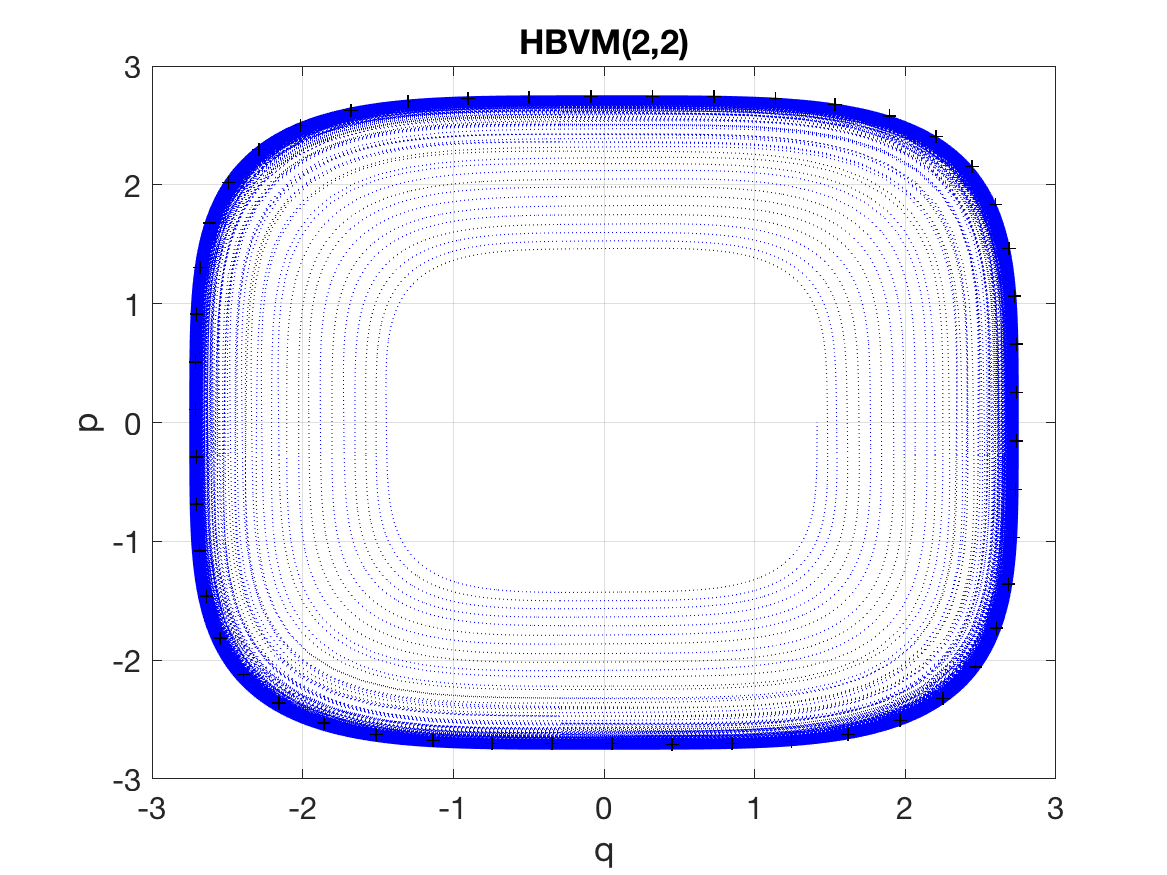}
\includegraphics[width=.45\textwidth]{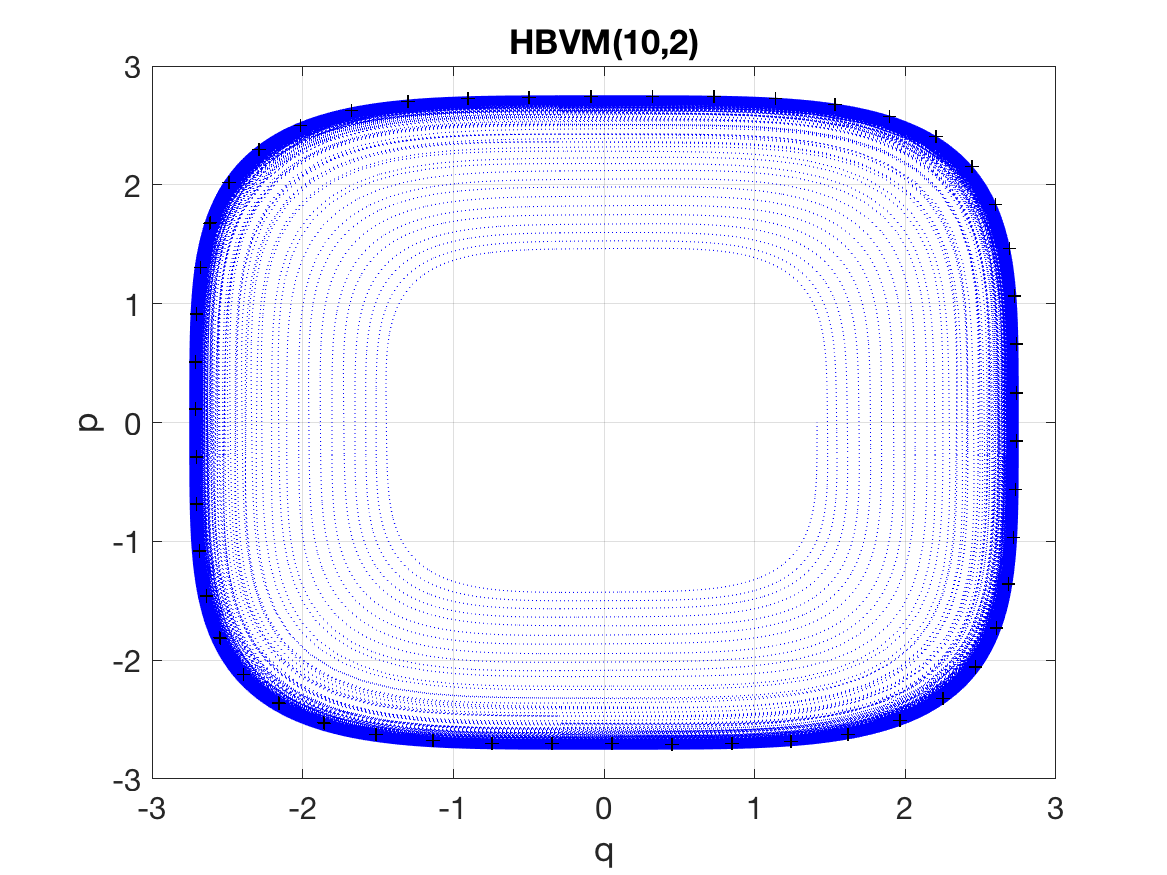}
\caption{Numerical results for problem \eqref{5.1}--\eqref{ICtest} with and \eqref{5.4} by using HBVM($2$,$2$) (left plots) and HBVM($10$,$2$) (right plots) in the time interval $[0,10^5]$, with step size $h=1/50$ (we refer to the text for more details).}
\label{fig1}
\end{center}
\end{figure}

\noindent
Additionally, Figure \ref{fig1} summarizes the results obtained by solving problem \eqref{5.4} with HBVM($2$,$2$) and HBVM($10$,$2$) in the time interval $\left[ 0, T={{\rm{10^5}}} \right]$, with step size $h=1/50$.  In the upper row of Figure \ref{fig1} are the plots of the numerical Hamiltonian $H(q_{n},p_{n})$ ($n=0,...,N=T/h$) with respect to time, from which one deduces that both methods fluctuate in the orders range $[10^{-8}, 10^{-2}]$, as more deeply inspected by the plots in the second row of the figure, where the computed values of $|\Delta H|=|H(q_n,p_n)-H(q_{n-1},p_{n-1})|$ ($n=1,2,...,N$) with respect to time are depicted. The third row of Figure \ref{fig1} then shows the numerical trajectories, in the phase space $(q,p)$, given by the two methods, in which the points within the time interval ${\left[ {{{10}^4} - 2,{{10}^4}} \right]}$ are marked by a black plus.

\noindent
As one may easily notice, Figure \ref{fig1} reveals analogous results for both the Gauss Legendre HBVM($2$,$2$) and the HBVM($10$,$2$); moreover, the same is true whether the corresponding plots via the spectrally-accurate higher order HBVM($22$,$20$) are considered, that we here omit for the sake of redundancy. Nevertheless, this is not always the case, as revealed by the numerical tests of the two following subsections, in which the advantage of recurring to HBVM($k$,$s$) methods with $k>s$ turns out to be evident.

\subsection{Problem 2}
Consider now \eqref{5.1}--\eqref{ICtest} with 

\begin{equation}\label{5.5}
	\begin{array}{lllll}
		m=1,~~~~~~\qquad H(q,p)=\frac{1}{2}(q^{2}-\cos q),\\
		\alpha=-10^{-5},~~~~~~q_{0}=0,~~~~ p_{0}=1.99999.	
	\end{array}
\end{equation}

\noindent
This problem can be seen as a dissipative delay-variant of the nonlinear pendulum, resulting into a delay Hamiltonian problem with dissipation of the Hamiltonian. In accordance with the analogous example in \cite{Brugnano22-1}, for the given initial conditions, the pendulum should undergo damped oscillations with a decreasing trend of the Hamiltonian function, being the numerical reproducing of the right dissipation of the Hamiltonian crucial, when relatively large step sizes are used. We preliminary solve the problem with HBVM($2$,$2$), HBVM($10$,$2$) and HBVM($15$,$3$), to verify the theoretical orders $4$ (fist two methods) and $6$ (latter one) of convergence. The results are reported in Table \ref{table2}, from which we highlight that, for each value of the step size $h$, HBVM($10$,$2$) achieves smaller errors $\varepsilon_N(h)$ than HBVM($2$,$2$).

\begin{table}[H]
     \centering
	
	\caption{Estimates of the last point error $\varepsilon_N(h)$, together with the convergence order $p$, of HBVM$(2,2)$, HBVM$(10,2)$ and HBVM($15$,$3$) to problem \eqref{5.1}--\eqref{ICtest} with \eqref{5.5}, referring to the time interval $[0,2]$ and step sizes $h= q^{-1}$ ($q = 2$, $4$, $8$, $16$, $32$).}
	\begin{tabular}{lllllll}
\hline
\multicolumn{1}{c}{\multirow{2}{*}{$h$}} & \multicolumn{2}{c}{HBVM($2$,$2$)}                                & \multicolumn{2}{c}{HBVM($10$,$2$)}                                & \multicolumn{2}{c}{HBVM($15$,$3$)}                             \\ \cline{2-7} 
\multicolumn{1}{c}{}                   & \multicolumn{1}{c}{$\varepsilon_N(h)$} & \multicolumn{1}{c}{$p$} & \multicolumn{1}{c}{$\varepsilon_N(h)$} & \multicolumn{1}{c}{$p$} & \multicolumn{1}{c}{$\epsilon_N(h)$} & \multicolumn{1}{c}{$p$} \\ \hline
1/2                                    & 1.3448e-04                             & \multicolumn{1}{c}{---} & 3.2262e-05                             & \multicolumn{1}{c}{---} & 4.8112e-07                          & \multicolumn{1}{c}{---} \\ \hline
1/4                                    & 6.9255e-06                             & 4.2794e+00              & 2.1365e-06                             & 3.9165e+00              & 7.7446e-09                          & 5.9571e+00              \\ \hline
1/8                                    & 4.1458e-07                             & 4.0622e+00              & 1.3545e-07                             & 3.9794e+00              & 1.2193e-10                          & 5.9890e+00              \\ \hline
1/16                                   & 2.5638e-08                             & 4.0153e+00              & 8.4961e-09                             & 3.9949e+00              & 1.9069e-12                          & 5.9987e+00              \\ \hline
1/32                                   & 1.5982e-09                             & 4.0038e+00              & 5.3148e-10                             & 3.9987e+00              & 2.9754e-14                          & 6.0020e+00              \\ \hline
\end{tabular}\label{table2}
\end{table}

\begin{figure}[H]
\begin{center}
\includegraphics[width=.45\textwidth]{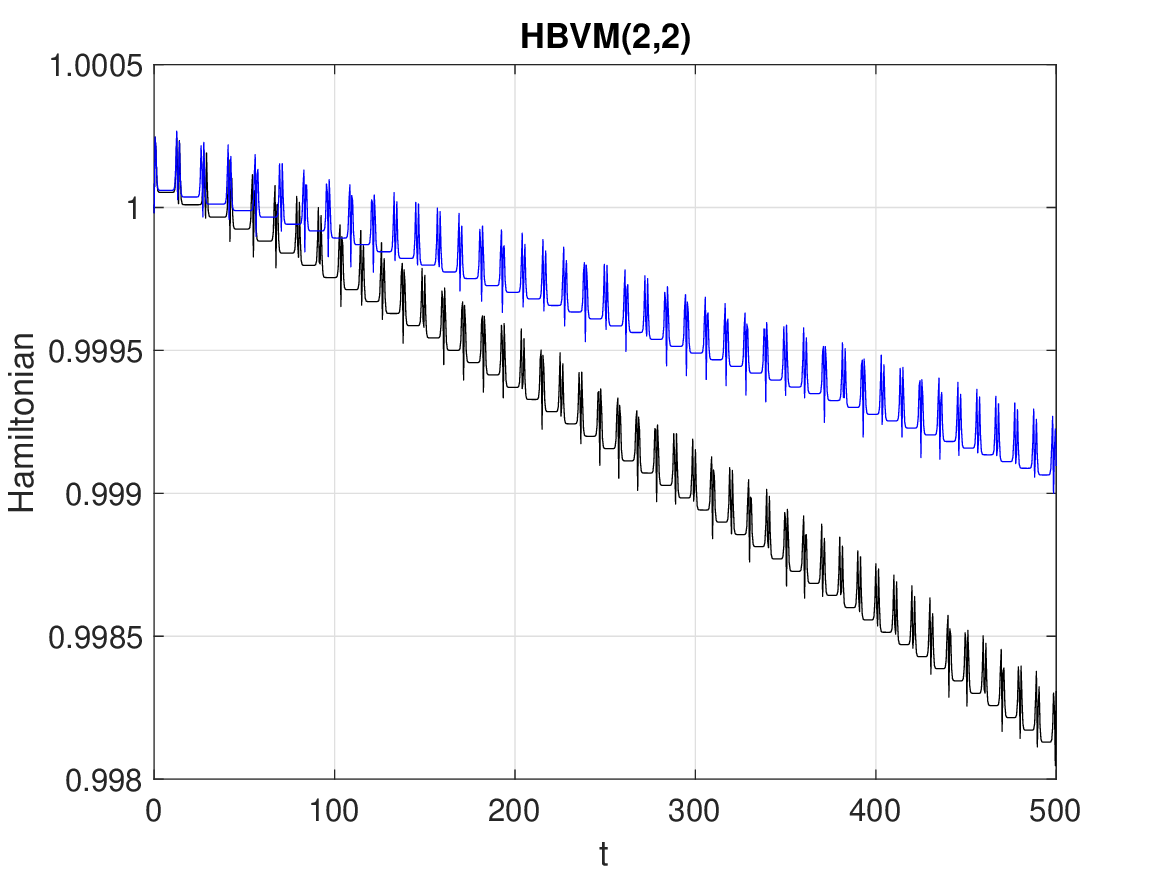}
\includegraphics[width=.45\textwidth]{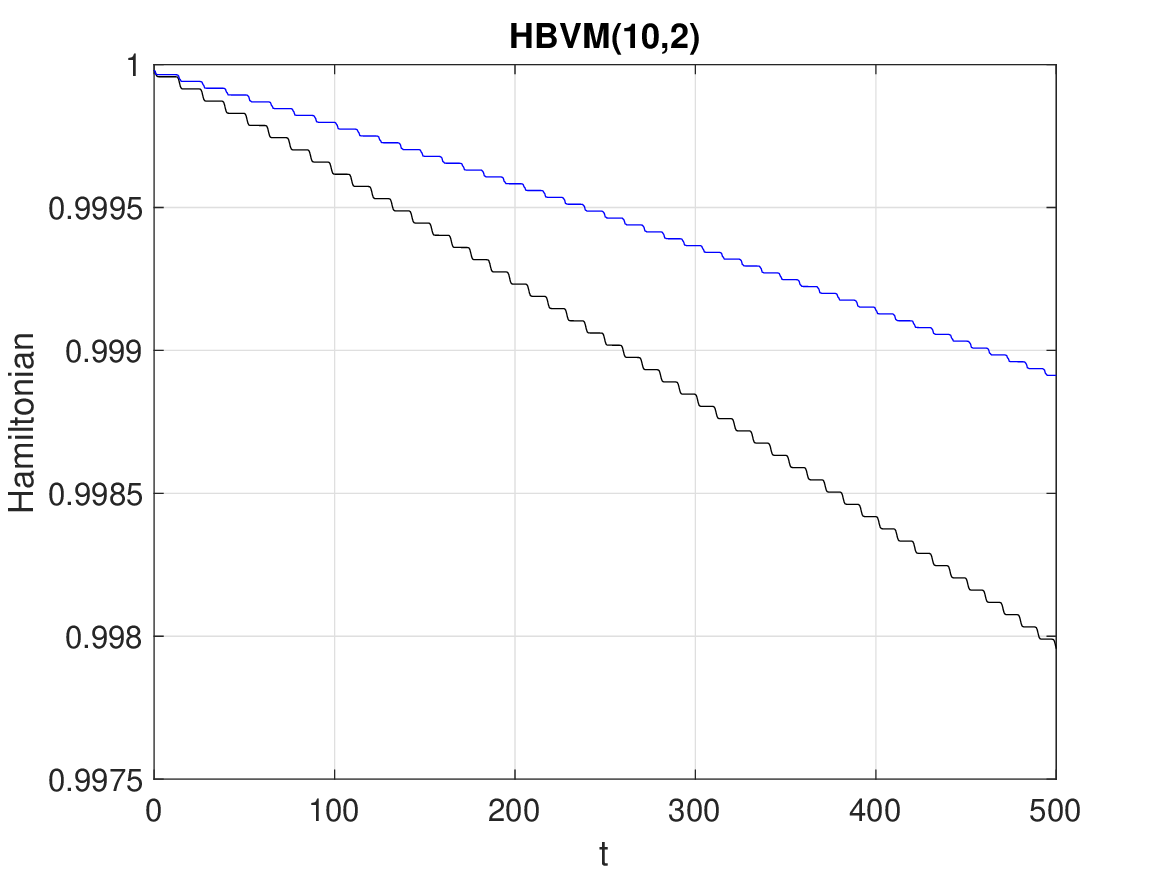}
\includegraphics[width=.45\textwidth]{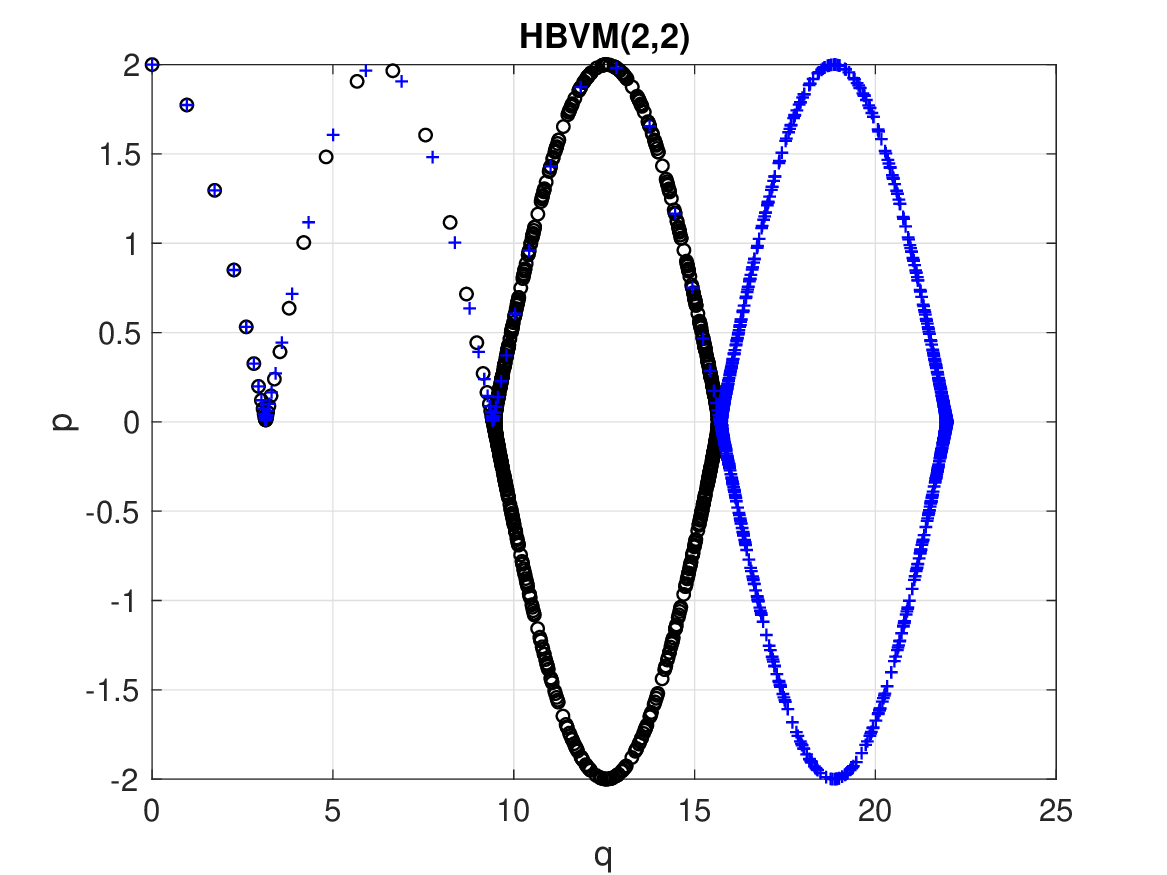}
\includegraphics[width=.45\textwidth]{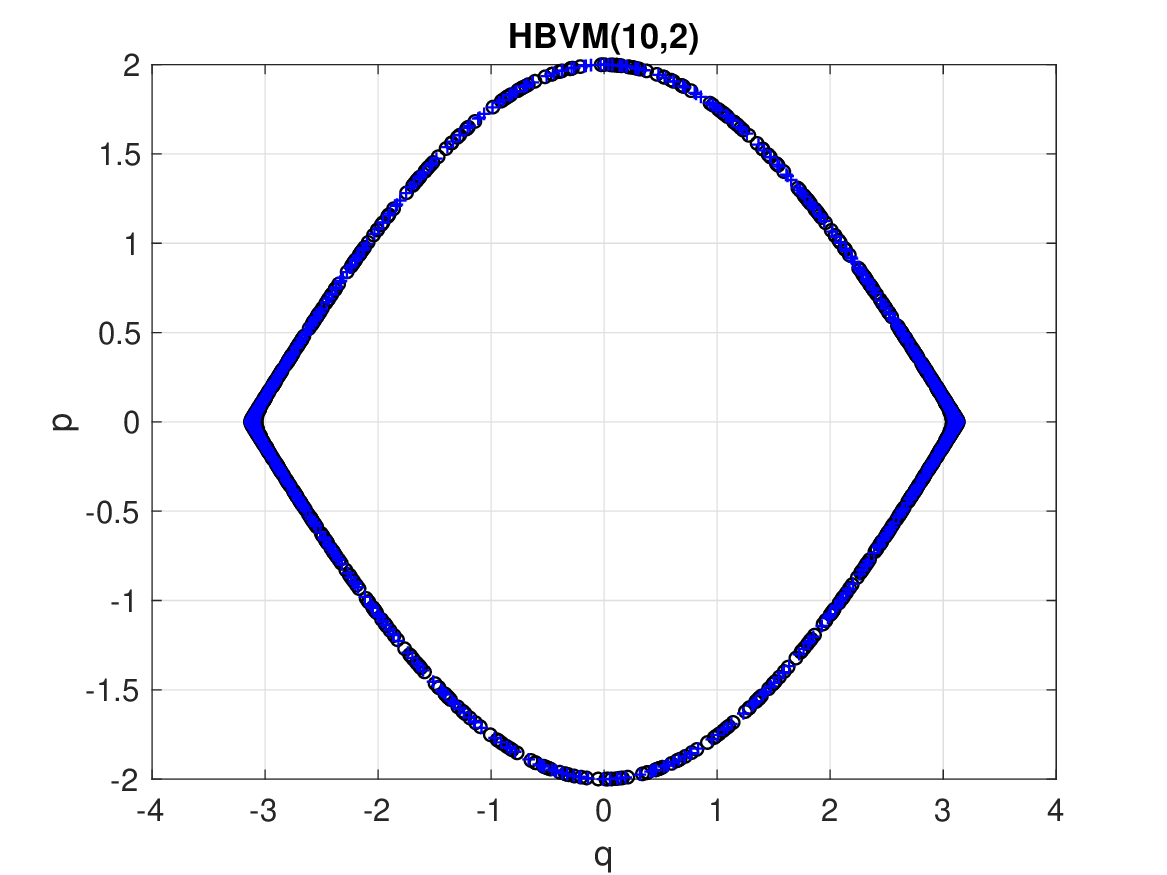}
\includegraphics[width=.45\textwidth]{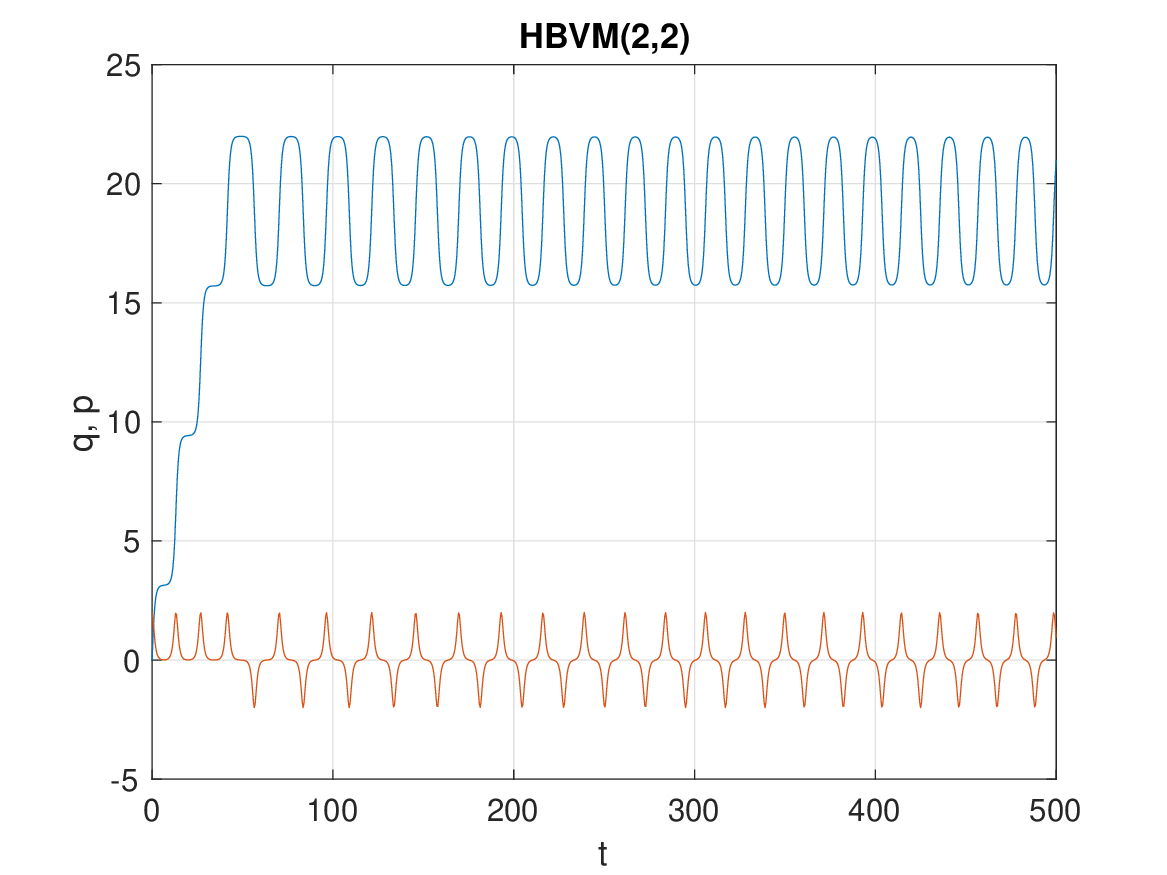}
\includegraphics[width=.45\textwidth]{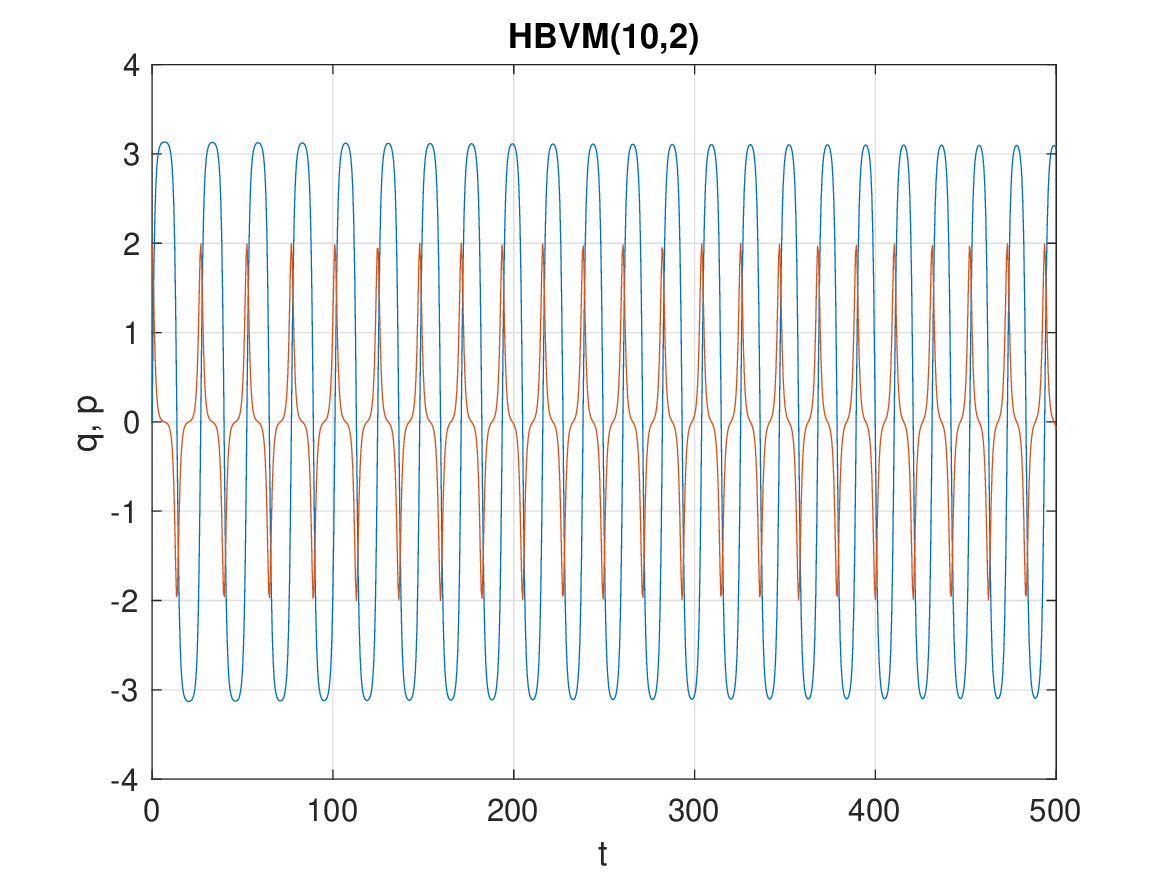}

	\caption{Numerical results for problem \eqref{5.1}--\eqref{ICtest} with \eqref{5.5} and \cite[Problem 3]{Brugnano22-1} by using HBVM($2$,$2$) (left plots) and HBVM($10$,$2$) (right plots) in the time interval $[0,500]$, with step size $h=1/2$ (we refer to the text for more details).}
	\label{fig2}
	\end{center}
\end{figure}

We then solve problem \eqref{5.5} on the time interval $[0,500]$, with time-step $h=1/2$, via HBVM$(2,2)$ and HBVM$(10,2)$. The obtained results are shown in Figure \ref{fig2}.

\begin{figure}[H]
\begin{center}
\includegraphics[width=.45\textwidth]{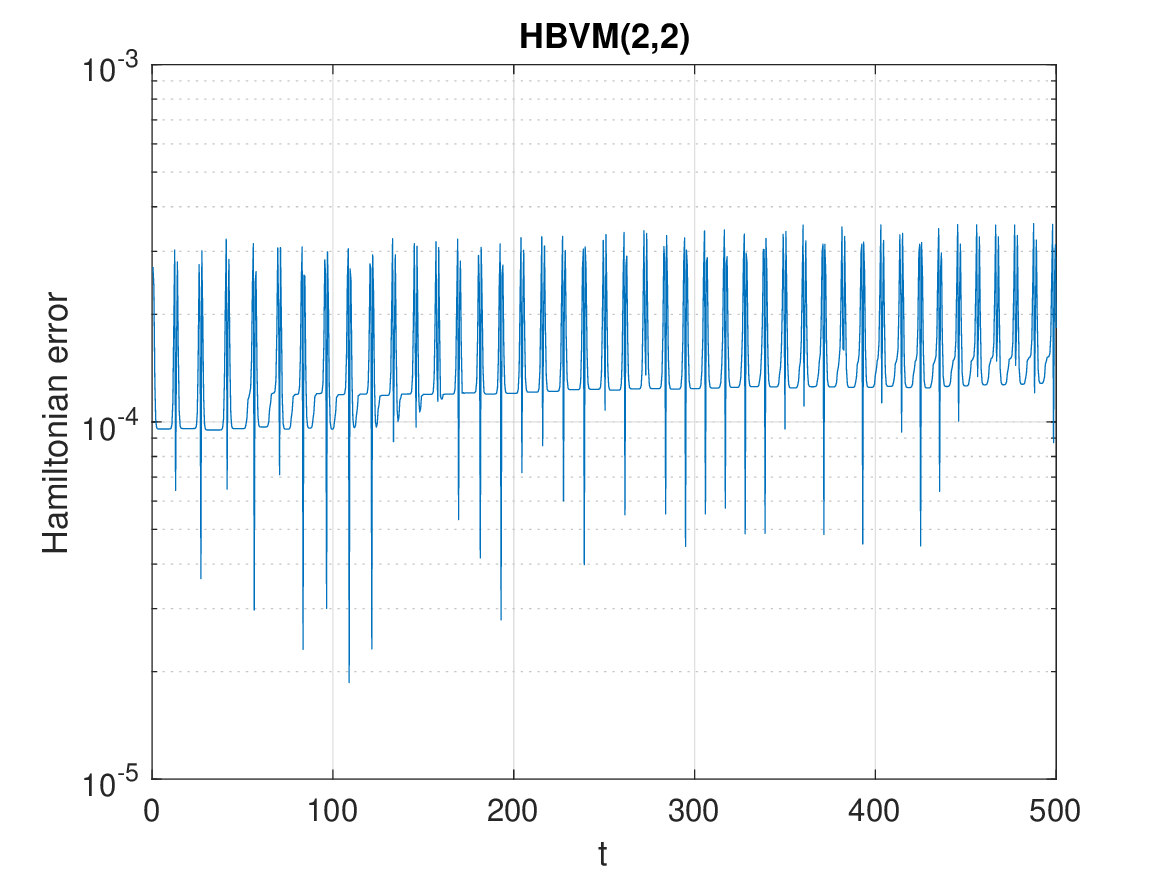}
\includegraphics[width=.45\textwidth]{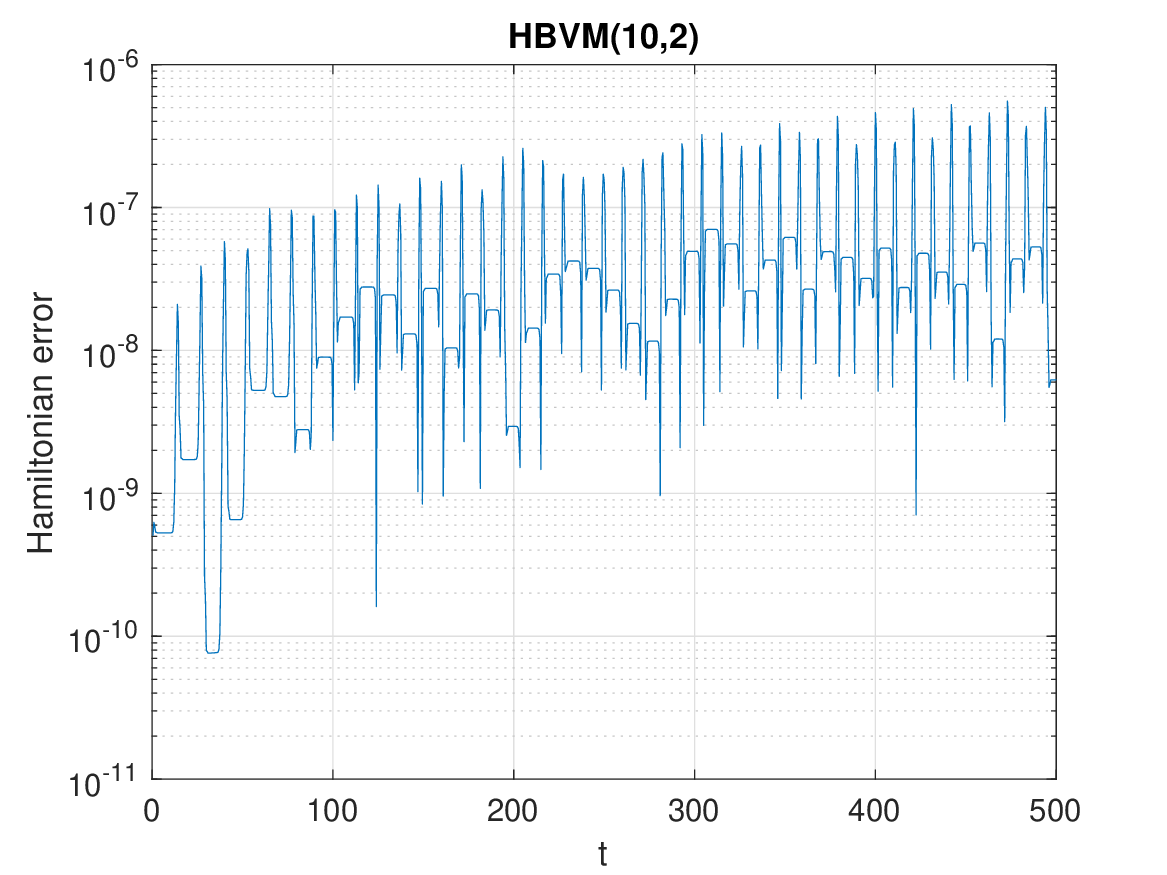}
\caption{Relative errors among time of the numerical Hamiltonian obtained by HBVM($2$,$2$) and HBVM($10$,$2$), when solving \eqref{5.1}--\eqref{ICtest} with \eqref{5.5},  with respect to the corresponding Hamiltonian values computed via HBVM($22$,$20$) (time interval $[0,500]$, step size $h=1/2$).}
\label{fig3}
\end{center}
\end{figure}

As one may easily see, the plots in the top row represent the numerical Hamiltonian $H(q_n,p_n)$ (blue line), revealing that both the methods exhibit a dissipation trend in the Hamiltonian function. However, for the HBVM$(2,2)$ method, the Hamiltonian values become significantly larger than $1$ for the smallest time instants, exhibiting fictitious oscillations. These oscillations cause the numerical solution to move away from the correct region of the phase space where the dynamics should occur. This is not the case of the resolution by HBVM$(10,2)$, that approximates the proper way of decreasing of the Hamiltonian, always remaining smaller than $1$. The more accurate computation of the Hamiltonian via HBVM($10$,$2$) is inspected in Figure \ref{fig3}, showing the plots of the relative errors on the Hamiltonian among time, computed with HBVM($2$,$2$) (i.e., the Gauss-Legendre method of order $4$) and HBVM($10$,$2$), with respect to the spectrally-accurate higher order HBVM($22$, $20$). As a deepening, the upper plots of Figure \ref{fig1} also report (black line) the corresponding Hamiltonian trend on the resolution of \eqref{5.1}--\eqref{ICtest} with \eqref{5.5} employing the continuous delay argument $(t-1)$ in place of the discontinuous $(\lfloor t\rfloor)$ in \eqref{5.1} and under the constant initial conditions $q(t)\equiv 0$ and $p(t)\equiv 1.9999$ (\cite[Problem 3]{Brugnano22-1}), obtained implementing the HBVM method of \cite{Brugnano22-1}. As expected, the smaller slope of the Hamiltonian decreasing trend of the FDEPCA case is outlined.

In the middle plots of Figure \ref{fig2}, the numerical solution of \eqref{5.1}--\eqref{ICtest} with \eqref{5.5} (blue points) is displayed in the phase space $(q,p$). It can be observed that the solution provided by HBVM$(2,2)$ ``jumps" three times before reaching the invariant region. This indicates that the pendulum completes three full rotations before losing enough energy to oscillate around its rest position. In contrast, the numerical solution obtained using the HBVM$(10,2)$ method always remains trapped in the correct region. For the sake of comparison with the above mentioned DDE case, the corresponding results related to \cite[Problem 3]{Brugnano22-1}, depicted by the black points, show analogous orbits, with the difference that the invariant region is reached after just two ``jumps'' of the trajectory by the HBVM($2$,$2$).

Finally, in the bottom plots of Figure \ref{fig2} is the numerical solution of \eqref{5.1}--\eqref{ICtest} with \eqref{5.5} with respect to time ($q(t)$ and $p(t)$ are represented by a blu and a red line, respectively), confirming that the HBVM$(2,2)$ method produces a solution that ``jumps" three times, while the solution obtained by the HBVM$(10,2)$ method does not exhibit such feature.

\subsection{Problem 3}
We finally face the following adapted Cassini ovals model (see, e.g., \cite{Brugnano16-1}), based on problem \eqref{5.1}--\eqref{ICtest} with:

	\begin{equation}\label{Problem3}
	     \begin{array}{l}
			m=1,~~H(q,p)=(q^{2}+p^{2})^{2}-10(q^{2}-p^{2}),~~\\
			 \alpha=10^{-5},~~~~q_{0}=0,~~~~ p_{0}=10^{-6}.	
		\end{array}
	\end{equation}

\noindent	
As expected, the theoretical order $2s$ of convergence of HBVM($k$,$s$) is numerically confirmed by the results in Table \ref{table3}, that takes into consideration HBVM($2$,$2$), HBVM($10$,$2$) and HBVM($15$,$3$).

\begin{table}[H]
     \centering
	
	\caption{Estimates of the last point error $\varepsilon_N(h)$, together with the convergence order $p$, of HBVM$(2,2)$, HBVM$(10,2)$ and HBVM($15$,$3$) to problem \eqref{5.1}--\eqref{ICtest} with \eqref{Problem3}, referring to the time interval $[0,35h]$ and step sizes $h=\frac{T}{100q}$ ($q = 1$, $2$, $4$, $8$, $16$).}
\begin{tabular}{lcccccc}
\hline
\multicolumn{1}{c}{\multirow{2}{*}{$h$}} & \multicolumn{2}{c}{HBVM($2$,$2$)} & \multicolumn{2}{c}{HBVM($10$,$2$)} & \multicolumn{2}{c}{HBVM($15$,$3$)} \\ \cline{2-7} 
\multicolumn{1}{c}{}                     & $\varepsilon_N(h)$  & $p$         & $\varepsilon_N(h)$   & $p$         & $\varepsilon_N(h)$   & $p$         \\ \hline
$T/100$                                  & 2.3787e-02          & ---         & 1.3856e-04           & ---         & 7.7321e-06           & ---         \\ \hline
$T/200$                                  & 1.3399e-05          & 1.0794e+01  & 1.1354e-05           & 3.6092e+00  & 1.6876e-07           & 5.5178e+00  \\ \hline
$T/400$                                  & 7.6790e-07          & 4.1251e+00  & 7.6292e-07           & 3.8955e+00  & 2.8672e-09           & 5.8791e+00  \\ \hline
$T/800$                                  & 4.8874e-08          & 3.9738e+00  & 4.8577e-08           & 3.9732e+00  & 4.5799e-11           & 5.9682e+00  \\ \hline
$T/1600$                                 & 3.0690e-09          & 3.9932e+00  & 3.0507e-09           & 3.9930e+00  & 7.2503e-13           & 5.9811e+00  \\ \hline
\end{tabular}\label{table3}
\end{table}

\begin{figure}[H]
\begin{center}
\includegraphics[width=.32\textwidth]{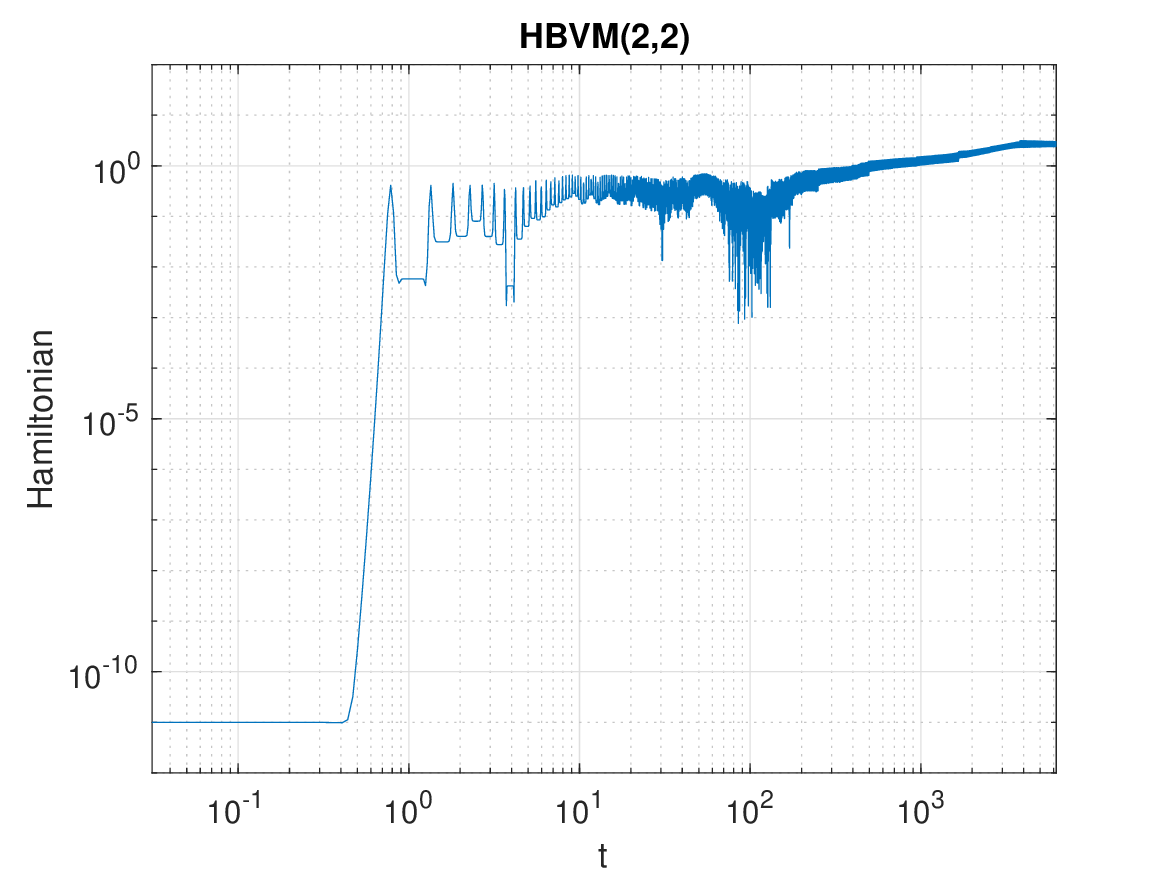}
\includegraphics[width=.32\textwidth]{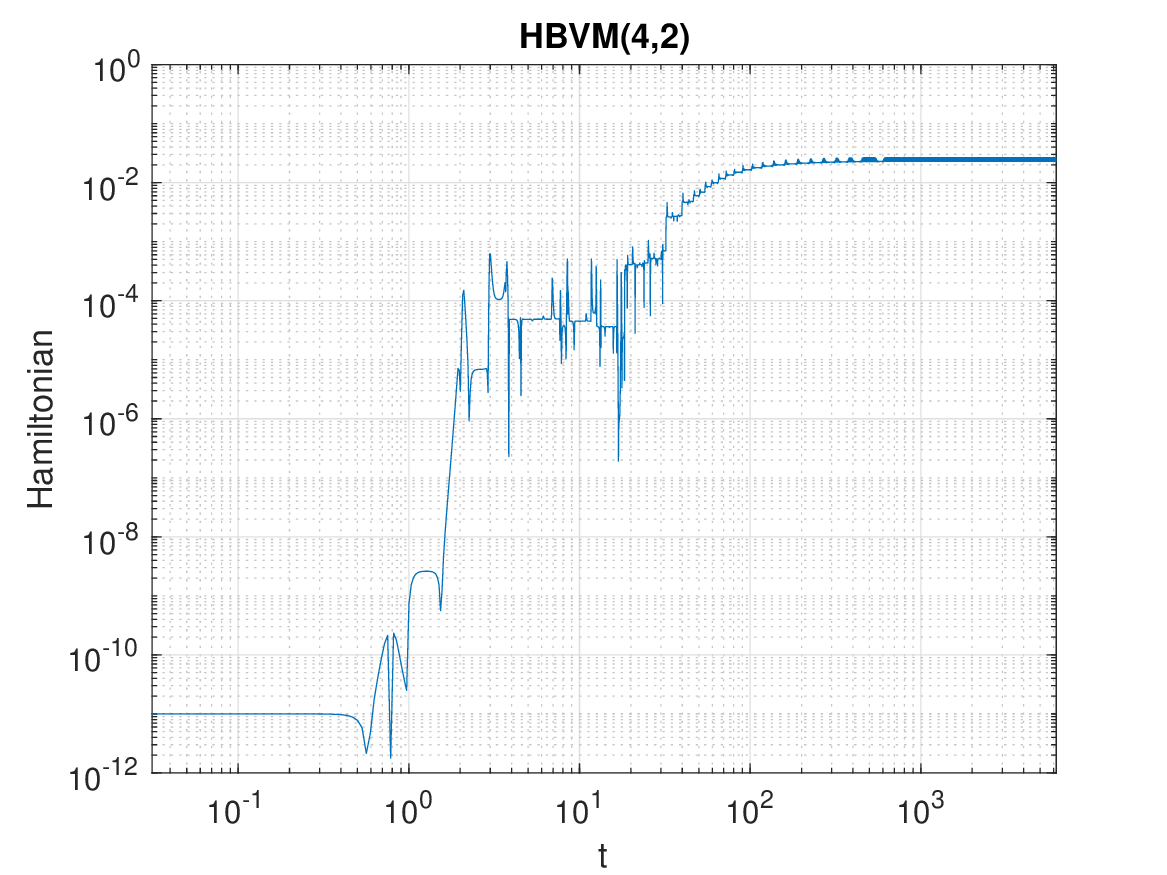}
\includegraphics[width=.32\textwidth]{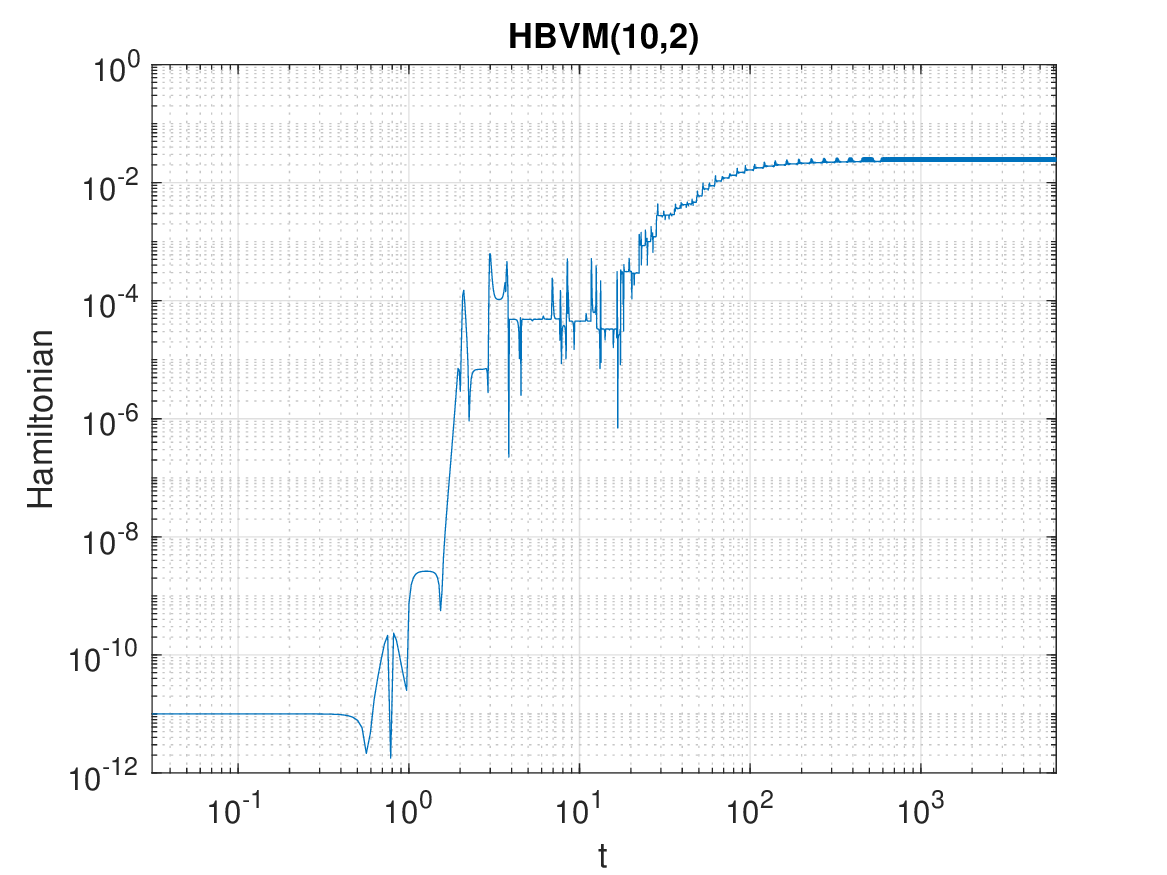}
\includegraphics[width=.32\textwidth]{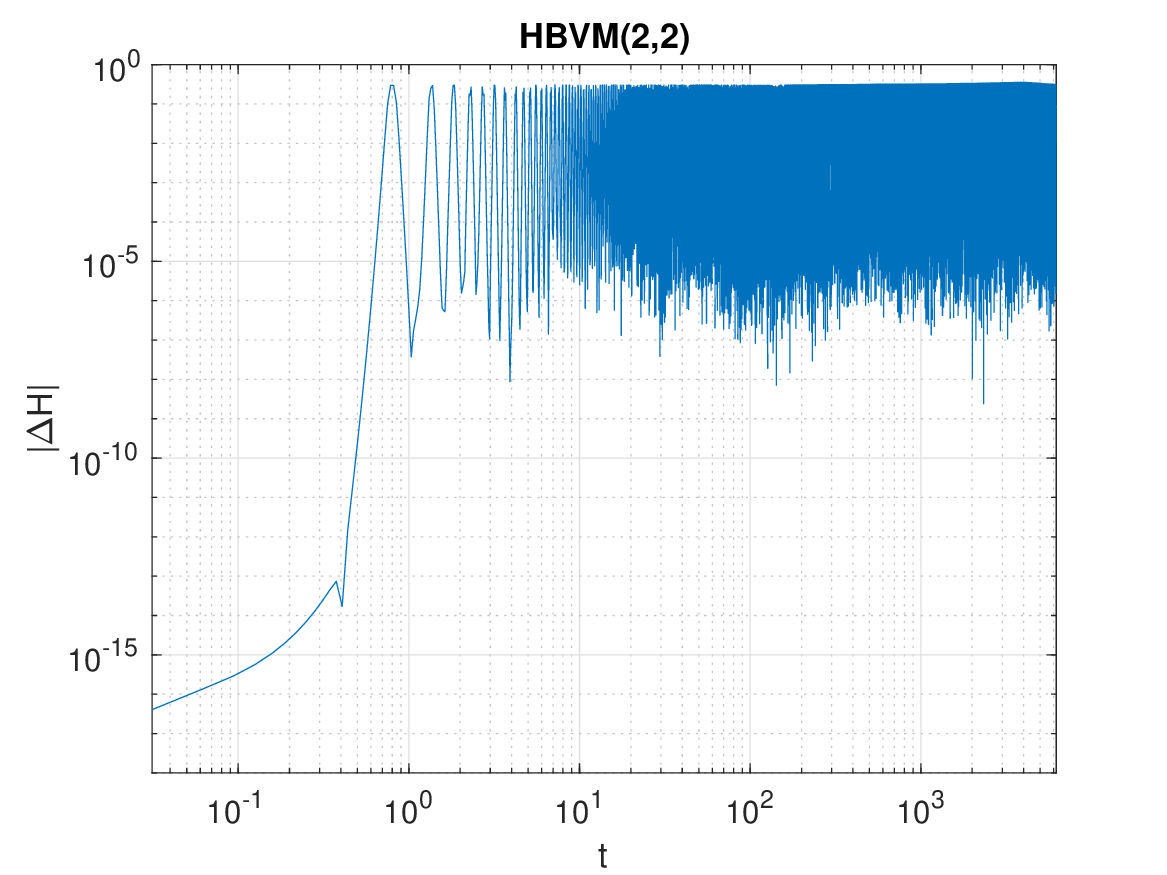}
\includegraphics[width=.32\textwidth]{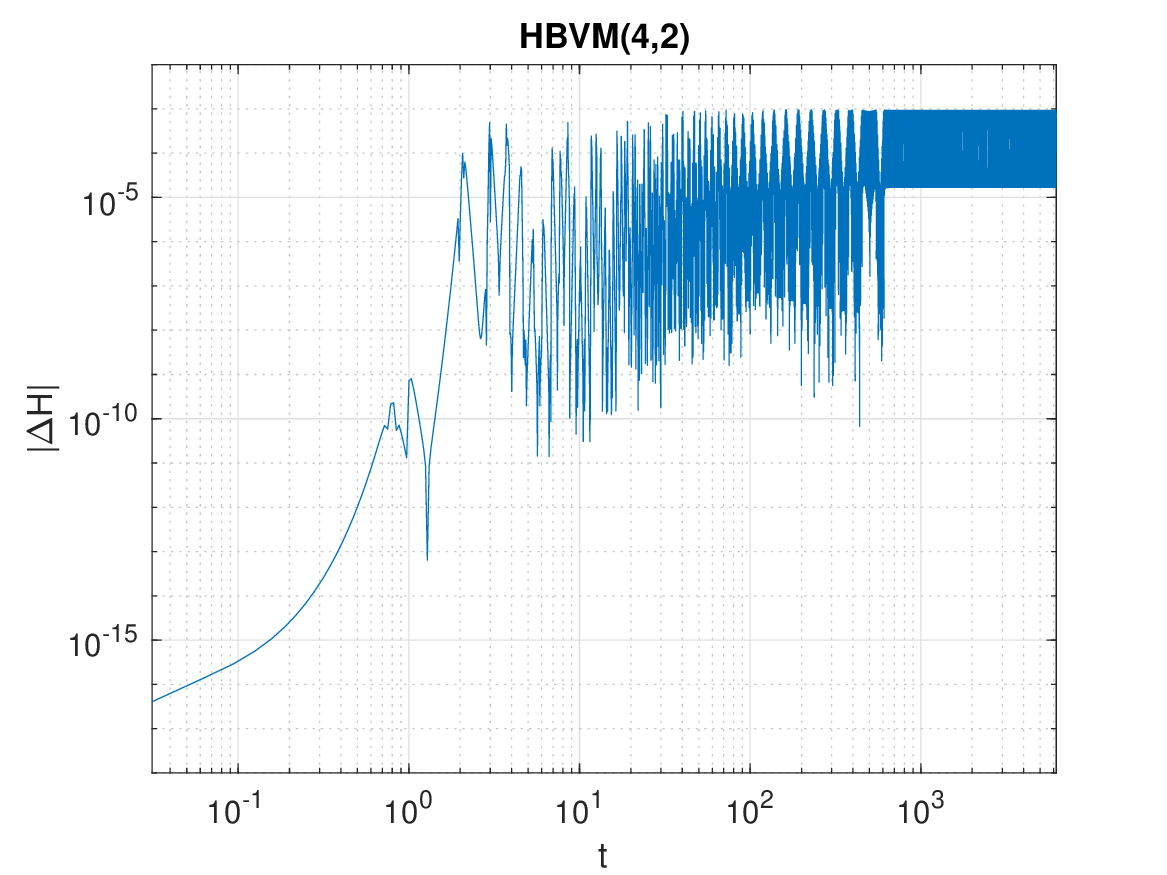}
\includegraphics[width=.32\textwidth]{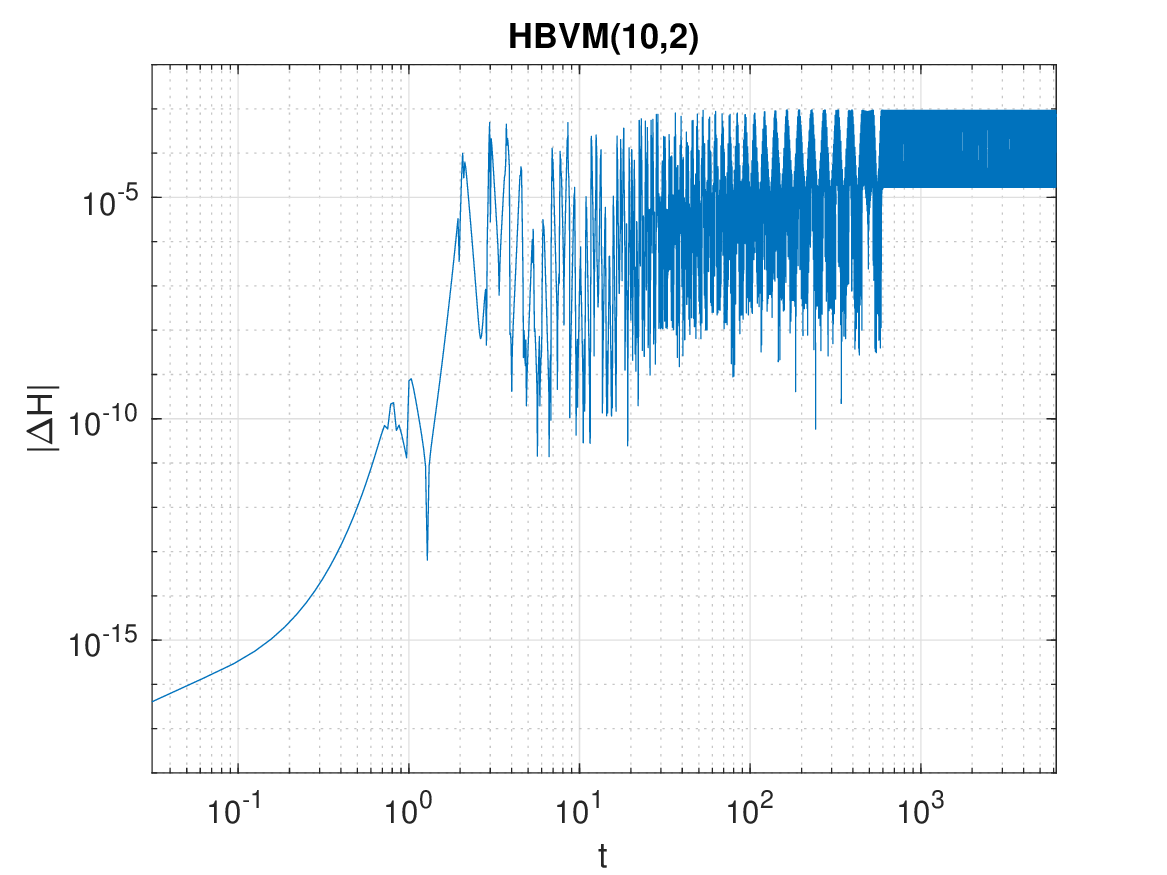}
\includegraphics[width=.32\textwidth]{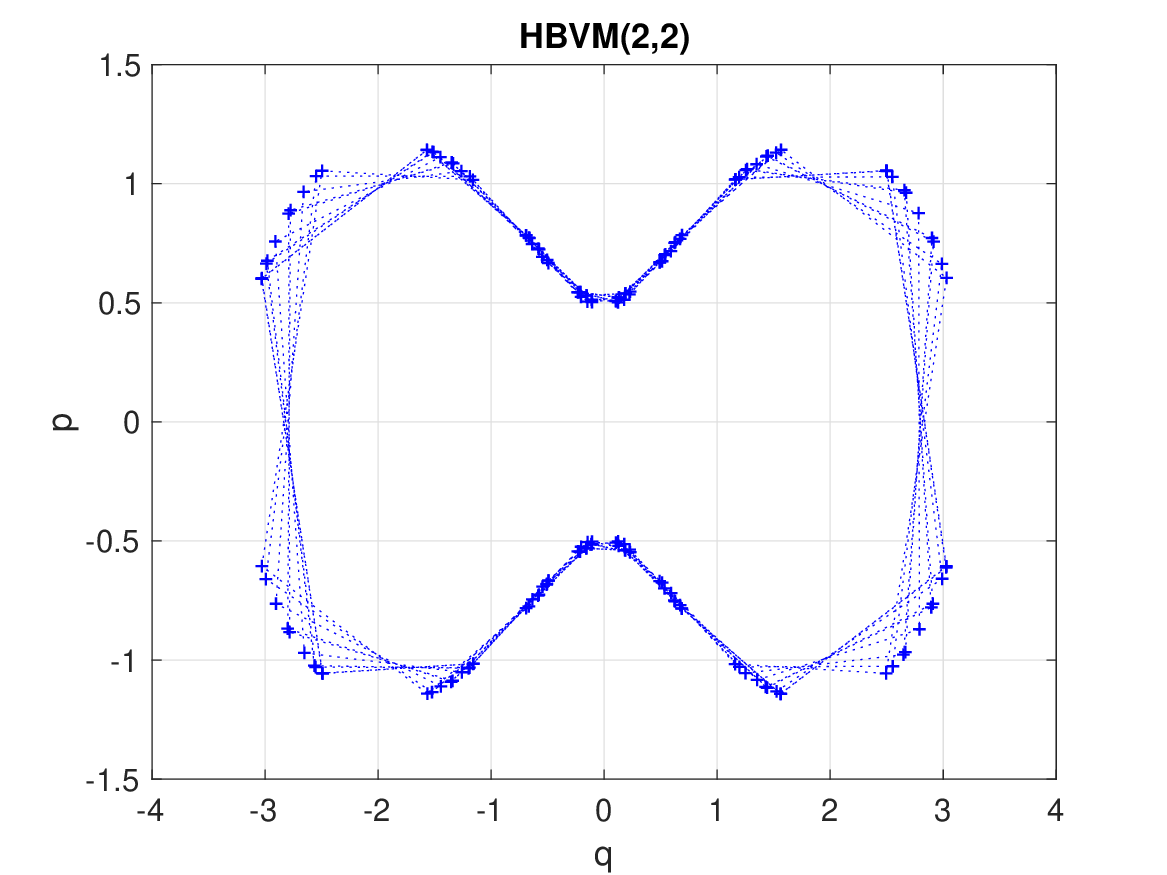}
\includegraphics[width=.32\textwidth]{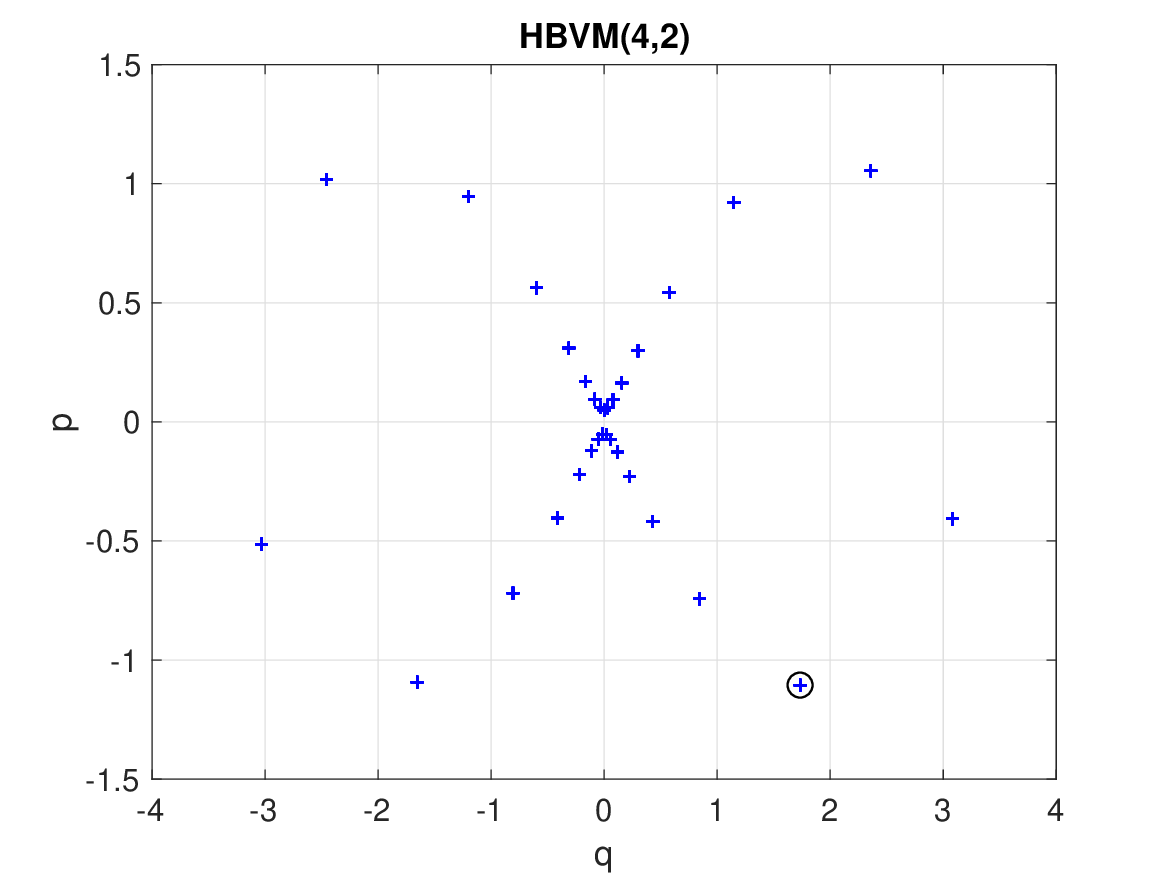}
\includegraphics[width=.32\textwidth]{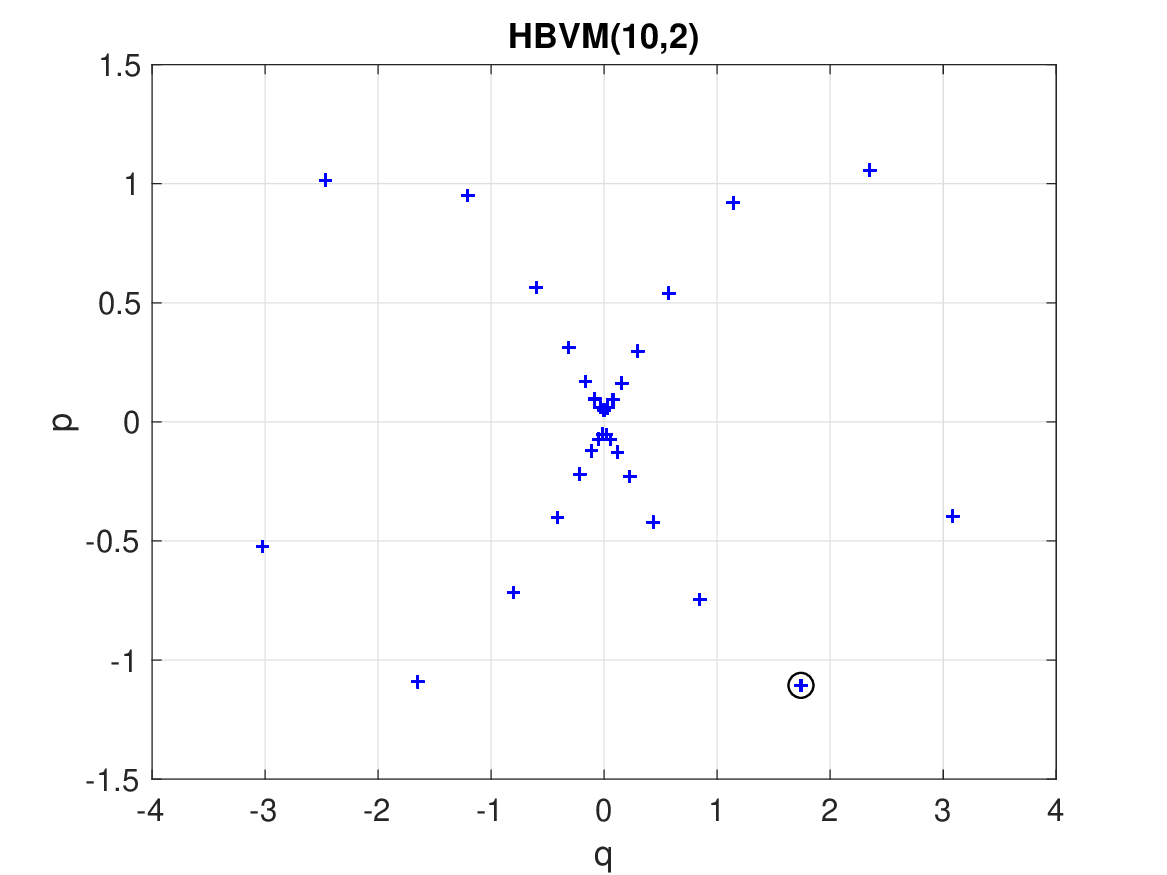}

	\caption{Numerical results for problem \eqref{5.1}--\eqref{ICtest} with \eqref{Problem3} by using HBVM($k$,$2$), for $k=2$ (left plots), $k=4$ (middle plots) and $k=10$ (right plots), in the time interval $[0,2\cdot 10^5h]$ (we refer to the text for more details).}
	\label{fig4}
	\end{center}
\end{figure}

\noindent
We solve problem \eqref{Problem3} by using HBVM$(k,2)$, for $k=2$, $4$, $10$, on the interval $[0,2\cdot 10^5h]$, with step size $h=T/100$, where $T=3.131990057003955$ is the period value in  \cite{Brugnano16-1}. The obtained results are shown in Figure \ref{fig4}.

In particular, in the upper row of Figure \ref{fig4} are the plots of the numerical Hamiltonian (in absolute value) $|H(q_{n},p_{n})|$ among time, from which one deduces that the HBVM($k$,$2$) with $k=4$ and $k=10$, quite soon reach a stationary behavior, while for $k=2$ (Gauss-Legendre) this is not the case.

In fact, by looking at the second row of plots in Figure \ref{fig4}, the asymptotic behavior of the three numerical solutions can be better discerned, revealing that only for HBVM($k$,$2$) with $k=4$ and $k=10$ the difference $|\Delta H|=|H(q_n,p_n)-H(q_{n-1},p_{n-1})|$ eventually reaches a significant smaller band of oscillation. Such a feature is completely analogous with what we get from the corresponding resolution via the spectrally-accurate HBVM($22$,$20$) method (whose pictures are omitted in Figure \ref{fig4} for the sake of redundancy).

The more likely reproduction of the Hamiltonian given by HBVM($4$,$2$) and HBVM($10$,$2$) then results into the classical Cassini oval periodic orbit in the phase space $(q,p)$, crossing the initial condition $(q_0,p_0)=(0,10^{-6})$, as reported in the two right-most plots at the bottom of Figure \ref{fig4}. In contrast, the orbit given by HBVM($2$,$2$), whose last $156$ points are depicted in the left-most plot at the bottom of Figure \ref{fig4}, is clearly not periodic. At this regard we emphasize that the choice a relatively small $k$ value larger than $s=2$ enables to gain the periodic feature of the orbit in the phase space. This is more deeply inspected in Table \ref{table4}, where the last $20$ points $(q,p)$ of the orbit computed after each period and lying inside the black circles in the plots at the bottom of Figure \ref{fig4} are displayed, thus confirming that both HBVM($4$,$2$) and HBVM($10$,$2$) are able to preserve at least the first $10$ significant digits.

\begin{table}[H]
  \centering
	\caption{The last $20$ points $(q,p)$ of the orbit, referring to problem \eqref{5.1}--\eqref{ICtest} with \eqref{Problem3}, computed after each period and lying inside the black circles in the plots at the bottom of Figure \ref{fig4}: HBVM($4$,$2$) (left-hand side data), HBVM($10$,$2$) (right-hand side data).}
\begin{tabular}{llll}
\hline
\multicolumn{2}{c}{HBVM($4$,$2$)}                                    & \multicolumn{2}{c}{HBVM($10$,$2$)}                                   \\ \hline
\multicolumn{1}{c}{$q$} & \multicolumn{1}{c}{$p$}                    & \multicolumn{1}{c}{$q$} & \multicolumn{1}{c}{$p$}                    \\ \hline
1.733257346741073  & \multicolumn{1}{c}{-1.105066941919119} & 1.741260676483831   & \multicolumn{1}{c}{-1.106100656753309} \\ \hline
1.733257346733711   & -1.105066941918150                     & 1.741260676498351   & -1.106100656755148
                     \\ \hline
1.733257346727258   & -1.105066941917299                     & 1.741260676513584   & -1.106100656757076
                     \\ \hline
1.733257346713944   & -1.105066941915546                     & 1.741260676527458   & -1.106100656758834
                     \\ \hline
1.733257346703106   & -1.105066941914119                     & 1.741260676538360   & -1.106100656760213
                     \\ \hline
1.733257346692424   & -1.105066941912712                     & 1.741260676550749   & -1.106100656761783
                     \\ \hline
1.733257346681209   & -1.105066941911235                     & 1.741260676560553   & -1.106100656763024
                     \\ \hline
1.733257346671835   & -1.105066941910001                     & 1.741260676568207   & -1.106100656763993
                     \\ \hline
1.733257346652067   & -1.105066941907396                     & 1.741260676576656   & -1.106100656765062
                     \\ \hline
1.733257346634478   & -1.105066941905079                     & 1.741260676578926   & -1.106100656765349
                     \\ \hline
1.733257346606794   & -1.105066941901432                     & 1.741260676578707   & -1.106100656765321
                     \\ \hline
1.733257346579902   & -1.105066941897892                     & 1.741260676588379   & -1.106100656766548
                     \\ \hline
1.733257346566756   & -1.105066941896163                     & 1.741260676595647   & -1.106100656767467
                     \\ \hline
1.733257346554184   & -1.105066941894506                     & 1.741260676590920   & -1.106100656766866
                     \\ \hline
1.733257346541494   & -1.105066941892836                     & 1.741260676580769   & -1.106100656765581
                     \\ \hline
1.733257346525058   & -1.105066941890670                     & 1.741260676575776   & -1.106100656764950
                     \\ \hline
1.733257346513068   & -1.105066941889092                     & 1.741260676577758   & -1.106100656765201
                     \\ \hline
1.733257346504652   & -1.105066941887984                     & 1.741260676589099   & -1.106100656766639
                     \\ \hline
1.733257346500618   & -1.105066941887454                     & 1.741260676602409   & -1.106100656768324
                     \\ \hline
1.733257346497192   & -1.105066941887002                     & 1.741260676613088   & -1.106100656769676
                     \\ \hline
\end{tabular}\label{table4}
\end{table}

\section{Conclusions and perspectives}\label{section7}
Moving from the derivation of the perturbation results related to the expansion of the vector field of the reference FDEPCA \eqref{Pb} along the Legendre orthonormal polynomial basis \eqref{Leg}, we have developed a class of high-order HBVM($k$,$s$) methods ($k\ge s\ge 1$) for \eqref{Pb}. The $2s$ order of the obtained family of methods is proved and its actual implementation is discussed. Numerical simulations for the resolution of Hamiltonian FDEPCA problems of a certain kind  have been performed, confirming the theoretical features of the methods. To the best of our knowledge, the present paper is the first one in which the HBVMs approach is applied to FDEPCAs. Moreover, the presented framework gives rise to generalizations along different directions. For one thing, HBVMs may be extended to a wider range of differential equations, such as integral delay differential equation, differential algebraic equation or fractional differential equation; in all these cases, the establishment of the corresponding perturbation theory deserves further to research and explore. For another,  finding approximations that belong to functional subspaces unlike polynomials could be of interest for future investigation. These two areas constitute worthwhile research ideas to be considered for future research.


\begin{thebibliography}{99}

\bibitem{Akhmet06} U. Akhmet, H. Oktem, W. Pickl, W. Weber, An anticipatory extension of malthusian model, AIP Conference Proceedings 839 (2006) 260.

\bibitem{Altintan06} D. Altintan, Extension of the logistic equation with piecewise constant arguments and population dynamics, Middle East Technical University, Turkey, 2006.

\bibitem{Brugnano 15-1} P. Amodio, L. Brugnano, F. Iavernaro, Energy-conserving methods for Hamiltonian boundary value problems and applications in astrodynamics, Adv.
Comput. Math. 41 (2015) 881-905.


\bibitem{Brugnano20-1}
 P. Amodio, L. Brugnano, F. Iavernaro, Analysis of Spectral Hamiltonian Boundary Value Methods
(SHBVMs) for the numerical solution of ODE problems, Numer. Algorithms. 83 (2020) 1489-1508.

\bibitem{Brugnano22-2} P. Amodio, L. Brugnano, F. Iavernaro, Continuous-Stage Runge-Kutta Approximation to differential
problems, 11 (2022) doc: https://doi.org/10.3390/axioms11050192.

\bibitem{Brugnano18-3} L. Barletti, L. Brugnano, G. Frasca Caccia, F. Iavernaro, Energy-conserving methods for the nonlinear Schr\"{o}dinger equation, Appl. Math. Comput. 318 (2018) 3-18.

\bibitem{BellZenn2003} A. Bellen, M. Zennaro, Numerical Methods for Delay Differential Equations, Clarendon Press, Oxford, 2003.

\bibitem{Bereketoglu10}H. Bereketoglu, G. Seyhan, A. Ogun, Advanced impulsive differential equations with piecewise constant arguments, Math. Model. Anal. 15 (2010) 175-187.

\bibitem{Brugnano97-1} L. Brugnano, Boundary value method for the numerical approximation of ordinary differential equa-
tions, Lect. Notes. Comput. Sc. 1196 (1997) 78-89.

\bibitem{Brugnano97-2} L. Brugnano, D. Trigiante, Block Boundary Value Methods for Linear Hamiltonian Systems, Appl.
Math. Comput. 81 (1997) 49-68.

\bibitem{Brugnano97-3} L. Brugnano, Essentially Symplectic Boundary Value Methods for Linear Hamiltonian Systems. J.
Comput. Math. 15 (1997) 233-252.

\bibitem{Brugnano98-1} L. Brugnano, D. Trigiante, Boundary value methods: the third way between linear multistep and Runge-Kutta methods, Math. Appl. Math. Comput. 36 (1998) 269-284.

\bibitem{Brugnano2002}
L. Brugnano, C. Magherini, Blended implementation of block implicit methods for ODEs, Applied numerical mathematics 42 (2002) 29-45.

\bibitem{Brugnano09-1} L. Brugnano, F. Iavernaro, D. Trigiante, Hamiltonian BVMs (HBVMs): A family of "drift-free" methods for integrating polynomial Hamiltonian systems. AIP Conference Proceedings 1168 (2009) 715-718.




\bibitem{Brugnano12-1} L. Brugnano, M. Calvo, J. I. Montijano, Energy preserving methods for Poisson systems.
J. Comput. Appl. Math. 236 (2012) 3890-3904.


\bibitem{Brugnano12-3} L. Brugnano, F. Iavernaro, Line integral methods which preserve all invariants of conservative problems, J. Comput. Appl. Math. 236 (2012) 3905-3919.

\bibitem{Brugnano12-4} L. Brugnano, F. Iavernaro, D. Trigiante, Energy and QUadratic invariants preserving integrators based upon Gauss collocation formulae, SIAM. J. Numer. Anal. 50 (2012) 2897-2916.

\bibitem{Brugnano14-1} L. Brugnano, Y. Sun, Multiple invariants conserving Runge-Kutta type methods for Hamiltonian
problems, Numer. Algorithms. 65 (2014) 611-632.



\bibitem{Brugnano16-1} L. Brugnano, F. Iavernaro, Line Integral Methods for Conservative Problems, CRC Press, 2016.

\bibitem{Brugnano18-1} L. Brugnano, G. Gurioli, F. Iavernaro, Analysis of energy and QUadratic invariant preserving (EQUIP) methods, J. Comput. Appl. Math. 335 (2018) 51-73.

\bibitem{Brugnano18-2} L. Brugnano, G. Gurioli, F. Iavernaro, E. Weinm\"{u}ller, Line integral solution of Hamiltonian systems with holonomic constraints, Appl. Numer. Math. 127 (2018) 56-77.

\bibitem{Brugnano18-4} L. Brugnano, F. Iavernaro, J. I. Montijano, L. R'andez, Spectrally Accurate Space-Time Solution of Hamiltonian PDEs, Numer. Algorithms. 81 (2019) 1183-1202.

\bibitem{Brugnano18-5} L. Brugnano, C. Zhang, D. Li, A class of energy-conserving Hamiltonian boundary value methods for nonlinear Schr\"{o}dinger equation with wave operator, Commun.  Nonlinear. Sci. Numer. Simul. 60 (2018) 33-49.

\bibitem{Brugnano18-6} L. Brugnano, F. Iavernaro, Line integral solution of differential problems, axioms. 7 (2018) 36. doc: https://doi.org/10.3390/axioms7020036.

\bibitem{Brugnano19-1} L. Brugnano, J. I. Montijano, L. R'andez, On the effectiveness of spectral methods for the numerical solution of multi-frequency highly-oscillatory Hamiltonian problems, Numer. Algorithms. 81 (2019) 345-376.


\bibitem{Brugnano19-3} L. Brugnano, G. Frasca Caccia, F. Iavernaro, Line integral solution of Hamiltonian PDEs, Mathematics-Basel 7 (2019) doc: https://doi.org/10.3390/math7030275.


\bibitem{Brugnano19-4} L. Brugnano, G. Gurioli, Y. Sun, Energy-conserving Hamiltonian boundary value methods for the numerical solution of the Korteweg-de Vries equation,
J. Comput. Appl. Math. 351 (2019) 117-135.

\bibitem{Brugnano19-5} L. Brugnano, G. Gurioli, C. Zhang, Spectrally accurate energy-preserving methods for the numerical solution of the "good" Boussinesq equation, Numer. Math. Part. D. E. 35 (2019) 1343-1362.

\bibitem{Brugnano22-1} L. Brugnano, G.Frasca-Caccia, F.Iavernaro, A new framework for polynomial approximation to differential equations, Adv. Comput. Math. 48 (2022) doc: https://doi.org/10.1007/s10444-022-09992-w.

\bibitem{Busenberg12} S. Busenberg, Cooke K, Vertically transmitted diseases: models and dynamics, Springer Science \& Business Media, 2012.

\bibitem{Chiu19}S. Chiu, T. Li, Oscillatory and periodic solutions of differential equations with piecewise constant generalized mixed arguments, Math. Nachr. 292 (2019) 2153-2164.

\bibitem{Cooke84} K. L. Cooke, J. Wiener, Retarded differential equations with piecewise constant delays, J. Math. Anal. Appl. 99 (1984) 265-297. doc:https://doi.org/10.1016/0022-247X(84)90248-8

\bibitem{Cooke86} K. L. Cooke, J. Wiener, Stability regions for linear equations with piecewise continuous delay, Comput. Math. Appl. 12 (1986) 695-701.

\bibitem{DosBar1986} J.G. Dos Reis, R. L. S. Baroni, On the existence of periodic solutions for autonomous retarded functional differential equations on $R^2$, Proceedings of the Royal Society of Edinburgh Section A: Mathematics 102.3--4 (1986) 259--262.


\bibitem{Esedog06} S. Esedog, Y-H.R. Tsai, Threshold dynamics for the piecewise constant Mumford-Shah functional, J. Comput. Phys. 211 (2006) 367-384.

\bibitem{Esmaeilzadeh20}M. Esmaeilzadeh, H. Najafi, H Aminikhah, A numerical scheme for diffusion-convection equation with piecewise constant arguments, comput. Methods. Differ. 8 (2020) 573-584.
       
\bibitem{Feng20}Z. Feng, Y. Wang, X. Ma, Asymptotically almost periodic solutions for certain differential equations with piecewise constant arguments, Adv. differ. equ-ny. (2020) doc:https://doi.org/10.1186/s13662-020-02699-6.

\bibitem{Gopalsamy92} K. Gopalsamy, Stability and oscillations in delay differential equations of population dynamics, Springer Science \& Business Media, 1992.

\bibitem{Gopalsamy98} K. Gopalsamy, P. Liu, Persistence and global stability in a population model, J. Math. Anal. Appl. 224 (1998) 59-80.

\bibitem{KapYor1974}
J.L. Kaplan, J.A. Yorke, Ordinary differential equations which yield periodic solutions of differential delay equations, J. Math. Anal Appl. 48 (1974) 317--324.

\bibitem{Kolmanovskii12} V. Kolmanovskii, A. Myshkis, Applied theory of functional differential equations, Springer Science \& Business Media, 2012.

\bibitem{Kumar13} P. Kumar, N. D. Pandey, D. Bahuguna, Existence of piecewise continuous mild solutions for impulsive functional differential equations with iterated deviating arguments, Electron. J. Differ. Eq. 241 (2013) 1-15.

\bibitem{Liang14} H. Liang, M. Liu, Z. Yang, Stability analysis of Runge-Kutta methods for systems $u'\left( t \right) = Lu\left( t \right) + Mu\left( {\lfloor t \rfloor} \right)$, Appl. Math. Comput. 228 (2014) 463-476.

\bibitem{Liu04} M. Liu, M. Song, Z. Yang, Stability of Rung-Kutta methods in the numerical solution of equation $u'\left( t \right) = au\left( t \right) + {a_0}u\left( {\lfloor t \rfloor} \right)$, J. Comput. Appl. Math. 166 (2004) 361-370.

\bibitem{Liu09} M. Liu, J. Gao, Z. Yang, Preservation of oscillations of the Runge-Kutta method for equation $x'\left( t \right) + ax\left( t \right) + {a_1}x\left(  \lfloor{t - 1}\rfloor  \right) = 0$, Comput. Math. Appl. 58 (2009) 1113-1125.

\bibitem{Liu15} X. Liu, M. Liu, Asymptotic stability of Runge-Kutta methods for nonlinear differential equations with piecewise continuous arguments, J. Comput. Appl. Math. 280 (2015) 265-274.

\bibitem{Li17} C. Li, C. Zhang, Block boundaryvalue methodsapplied tofunctional differentialequations withpiece-
wise continuous argument, Appl. Numer. Math. 115 (2017) 214-224.

\bibitem{Li19}C. Zhang, C. Li, Y. Jiang, Extended block boundary value methods for neutral equations with piece-
wise constant argument, Appl. Numer. Math. 150 (2019) 182-193.

\bibitem{Lv11} W. Lv, Z. Yang, M. Liu, Numerical stability analysis of differential equations with piecewise constant arguments with complex coefficients, Appl. Math. Comput. 218 (2011) 45-54.

\bibitem{MallNuss2011}
J. Mallet-Paret, R.D. Nussbaum, Stability of periodic solutions of state-dependent delay-differential equations, J Diff. Equ. 250 (2011) 4085-4103.

\bibitem{Mysh1977} A.D. Myshkis, On certain problems in the theory of differential equations with deviating arguments, Uspekhi Mat. Nuuk 32 (1977) 173-202.

\bibitem{Nuss1973}
R.D. Nussbaum, Periodic solutions of some nonlinear, autonomous functional differential equations, Bull. Amer. Math. Soc. 79 (1973) 811-814.

\bibitem{Nuss1979}
R.D. Nussbaum, Uniqueness and nonuniqueness for periodic solutions of $x'(t)=-g(x(t-1))$. J Diff. Equ. 34 (1979) 25--54.

\bibitem{Driver63} D. Rodney, A functional-differential system of neutral type arising in a two-body problem of classical electrodynamics, Academic Press, 1963.

\bibitem{Sobolev64} L. Sobolev, Partial differential equations of mathematical physics, Courier Corporation, 1964.

\bibitem{Song10} M. Song, X. Liu, The improved linear multistep methods for differential equations with piecewise continuous arguments, Appl. Math. Comput. 217 (2010) 4002-4009.

\bibitem{Wal1975}
H.-O. Walther, Existence of a non-constant periodic solution of a nonlinear autonomous functional
differential equation representing the growth of a single species population, J. Math. Biol. 1 (1975) 227--240.

\bibitem{Wang08} W. Wang, S. Li, Dissipativity of Runge-Kutta methods for neutral delay differential equations with piecewise constant delay, Appl. Math. Lett. 21 (2008) 983-991.

\bibitem{Wang13} W. Wang, Stability of solutions of nonlinear neutral differential equations with piecewise constant delay and their discretizations, Appl. Math. Comput.
219 (2013) 4590-4600.

\bibitem{Wen05} L. Wen, S. Li, Stability of theoretical solution and numerical solution of nonlinear differential equations with piecewise delays, J. Comput. Math. 23 (2005) 393-400.

\bibitem{Wiener93} J. Wiener, Generalized Solutions of Differential Equations, World Scientific, Singapore, 1993.

\bibitem{Zhang21} C. Zhang, X. Yan Convergence and stability of extended BBVMs for nonlinear delay-differential-algebraic equations with piecewise continuous arguments, Numer. Algorithms. 87 (2021) 921-937.

\end{thebibliography}
\end{document}